\documentclass[a4paper,twoside,12pt]{article}
\usepackage[english]{babel}
\usepackage{graphicx}
\usepackage{bbding}
\usepackage[latin1]{inputenc}
\usepackage{kpfonts}
\usepackage{fancyhdr}
\usepackage{amsmath,amssymb,color,epsfig}
\usepackage{mathrsfs}
\usepackage{eufrak}
\usepackage{dsfont}
\usepackage{upgreek}
\usepackage{fancyhdr}
\usepackage[dvipsnames]{xcolor}
\newtheorem{theorem}{Theorem}[section]
\newtheorem{proposition}[theorem]{Proposition}

\newtheorem{remark}[theorem]{Remark}
\newtheorem{lemma}[theorem]{Lemma}

\def\tr{{\rm tr}}

\def\1{\mathds{1}}
\def\e{\varepsilon}

\title{Density of Positive Eigenvalues of the Generalized Gaussian Unitary Ensemble }
\date{}
\definecolor{mon}{rgb}{10,126,140}
\begin{document}
\author{Mohamed BOUALI}
\maketitle
\begin{abstract}
We compute exact asymptotic of the statistical density  of random matrices belonging to the Generalized Gaussian orthogonal, unitary
and symplectic ensembles such that there no eigenvalues in the interval $[\sigma, +\infty[$. In particular, we show that the probability that all the eigenvalues of an
$(n\times n)$ random matrix are positive (negative) decreases for large $n$ as $\sim exp[-\beta\theta(\alpha)n^2]$ where
the Dyson index $\beta$ characterizes the ensemble, $\alpha$ is some extra parameter and the exponent $\theta(\alpha)$ is a function of $\alpha$ which will be given explicitly. For $\alpha=0$, $\theta(0)= (\log 3)/4 = 0.274653...$ is
universal. We compute the probability that the eigenvalues lie in the interval $[\sigma,+\infty[$ with $(\sigma>0,\; {\rm if }\;\alpha>0)$ and $(\sigma\in\mathbb R,\; {\rm if }\;\alpha=0)$. This generalizes the celebrated Wigner semicircle
law to these restricted ensembles. It is found that the density of eigenvalues generically exhibits
an inverse square-root singularity at the location of the barriers. These results generalized the case of Gaussian random matrices ensemble studied in \cite{D}, \cite{S}.\\
{\bf Math Subject classification:} 15B52, 15B57, 60B10.\\
{\bf Key-words:} Random matrices, Probability measures, Logarithmic potential.
\end{abstract}

\section{Introduction}
Random matrix theory has been successfully applied in various branches of physics and mathematics,
including in subjects ranging from nuclear physics, quantum chaos, disordered systems, and
number theory. Of particular importance are Generalized Gaussian random matrices with density is a Gaussian function times a power of the determinant. there are three classes of matrices distributed with such density: $(n \times n)$ real symmetric (Generalized Gaussian Orthogonal Ensemble (GGOE)), $(n \times n)$ complex Hermitian (Generalized Gaussian
Unitary Ensemble (GGUE)) and $(2n \times 2n)$ self-dual Hermitian matrices (Generalized Gaussian Symplectic Ensemble (GGSE)). In
these models the probability distribution for a matrix $X$ in the ensemble is given by
$$p_n(X)\propto |\det(X)|^{\beta\mu}\exp\Big(-\frac\beta 2\big< X,X\big>\Big).$$
where $\big<X,X\big>$ is the inner product on the space of matrices invariant, under orthogonal, unitary and symplectic
transformations respectively and the parameter $\beta$ is the Dyson index. In these three cases the inner products is defined as follow:
$$\begin{aligned}
&\big<X,X\big>=\tr(X^2),\;\beta=1\\
&\big<X,X\big>=\tr(X^*X),\;\beta=2\\
&\big<X,X\big>=\tr(X^\star X),\;\beta=4.\\
\end{aligned}$$
where $.^*$ denotes the hermitian conjugate of complex valued matrices and  $.^\star$ denotes the symplectic conjugate on
quaternion valued matrices.
A crucial result in the theory of random matrices is the celebrated Wigner semi-circle law. It states that for large
$n$ and on an average, the $n$ eigenvalues rescaled by the factor $\frac1{\sqrt n}$, lie within a finite interval $[-\sqrt 2, \sqrt 2]$, often referred to as the Wigner 'sea'. Within this sea, the statistical density of eigenvalues has a semi-circular form that vanishes at the two edges $-\sqrt 2$, $\sqrt 2$.
$$\rho(\lambda)=\frac1\pi\sqrt{2-\lambda^2}.$$
The above result means that, if one looks at the statistical density of eigenvalues of a typical system described by one of the three
ensembles above, for a large enough $n$, it will resemble closely to the Wigner semi-circle law.
From the semi-circle law, we know that on an average half the eigenvalues are positive and half of them are negative.

One of the main questions posed in \cite{A} and \cite{S} is: for a $n\times n$ Gaussian random matrix, what is the probability $P_n$ that
all its eigenvalues are positive (or negative).
 $$P_n={\rm Prob}[\lambda_1\geq 0,...,\lambda_n\geq 0].$$
 This probability has also been studied in the mathematics literature \cite{De} and one can easily compute $P_n$ for smaller values, $n=1, 2, 3$
$$P_1=\frac 12,\quad P_2=\frac12-\frac{\sqrt 2}{4},\quad P_3=\frac14-\frac{\sqrt 2}{2\pi}.$$
The interesting question is how $P_n$ behaves for large $n$? It was argued in \cite{A} that for large $n$, $P_n$ decays as $P_n\sim\exp\big(-\theta(0)n^2\big)$ where the decay constant $\theta(0)$ was estimated to be$\approx\frac14$ numerically and via a heuristic
argument. In \cite{S} the authors shown that for all the three Gaussian ensembles, to leading order in large $n$
$$P_n\sim\exp\big(-\beta\theta(0)n^2\big),\;\; {\rm where}\;\theta(0)=\frac{\log 3}{4}=0.2764....$$
 In our case we consider the same question in more general setting, where the ensemble of matrices is the Generalized Gaussian unitary ensemble equipped with the density $p_n(X)$ as in the previous. We will prove that, to leading order in large $n$, and for $\mu:=\mu_n\approx\alpha n$,  $(\alpha\geq 0)$
 $$P_n\sim\exp\big(-\beta\theta(\alpha)n^2\big),$$
 where $\theta(\alpha)$ is a function of $\alpha$ which will be given explicitly in proposition 3.5. In the simplest case for $\alpha\approx 0$ $(\alpha\leq 0.34)$,
$$\theta(\alpha)=\beta\big(\frac{\log3}{4} - C\alpha\big)+o(\alpha),$$
 where $\displaystyle C=\frac12+\frac{\log 2}2-\frac1{864}\Big(-36 (-6 + \sqrt6) + (54 - 161 \sqrt6) \log2 +
    27 (10 + \sqrt6) \log3\Big)\approx 0.3482$ and $o(\alpha)$ is a small terms in $\alpha$.

   Another important question will be studied in this work, namely  what is the statistical density of the negative (or positive) eigenvalues.

   In this paper we will calculate the asymptotic density of eigenvalues in this conditioned ensemble and we will see that it is quite different to the Wigner semi-circle law. We will prove the following result.
    For a positive real sequence $\mu_n$, such that $\mu_n\approx n\alpha$, $(\alpha\geq 0)$ and after scaling the statistical density of eigenvalues by $\frac1{\sqrt n}$, it converge as $n$ goes to $+\infty$ to the density $f_{\alpha, a}$ with support $[a,b]$, such that: for $\alpha>0$  $$f_{\alpha, a}(x)=\frac{1}{2\pi}\sqrt{\frac{b-x}{x-a}}\Big(2x+b-a-2\alpha\sqrt{\frac{a}{b}}\frac1x\Big),$$
   where $a>0$ and $b>a$ are the unique solutions of the following equations
   $$b+a-\frac{2\alpha}{\sqrt {ab}}=0,\;\;\frac{3}{4}(b-a)^2+a(b-a)+2\alpha\frac{\sqrt a}{\sqrt b}-2\alpha-2=0.$$
   It will be proved that for given $\alpha\geq 0$, the previous equations has a unique solutions $a_c>0$ and $\displaystyle b_c=\frac23\Big(\sqrt{6(\alpha+1)-2a_c^2}-\frac{a_c}2\Big)$.

   For $\alpha=0$
   $$f_{0, 0}(x)=\frac{1}{2\pi}\sqrt{\frac{b-x}{x}}\big(2x+b\big),$$
   where $a=0$ and $\displaystyle b=\frac23\sqrt 6.$

   More general for $\alpha\geq 0$ and $\sigma>0$. We will prove that the statistical density of eigenvalues such that there are no eigenvalues in the interval $]-\infty,\sigma]$ scaled by $\frac1{\sqrt n}$ converge to some density $f_{\alpha, a}$ with support $[a,b]$, where $a>0$ and $b>a$ are the unique solution in the interval $[\sigma,+\infty[$ of the following equations
   \begin{equation}\label{231} b+a-\frac{2\alpha}{\sqrt {ab}}\geq 0,\;\;\frac{3}{4}(b-a)^2+a(b-a)+2\alpha\frac{\sqrt a}{\sqrt b}-2\alpha-2=0.\end{equation}
   For $\alpha=0$ it is assumed that $\sigma\in\Bbb R$.

    More precisely for $\alpha >0$, if $\sigma\in[0, a_c]$, then $a=a_c$ and $b=b_c$. If $\sigma>a_c$, then $a=\sigma$ and $b$ is the unique solution of (\ref{231}).

    If $\alpha=0$, and $\sigma> -\sqrt 2$, then $a=\sigma$ and $b=\frac23\big(\sqrt{a^2+6}+\frac a2\big)$, if $\sigma\leq -\sqrt 2$, then $a=-\sqrt 2$ and $b=\sqrt 2$.

    For example for $\sigma\leq -\sqrt 2$, one recover's the famous Wigner semi-circle law.
    One can se that for $\alpha>0$, $\sigma$ must be strictly positif. This can be explain because the singularity at $0$ and there no single support of the density of eigenvalues when $\alpha >0$ and $\sigma<0$. This latter case will be studied by the author in forthcoming paper.

   The paper is organized as follows. In the first section we begin by recalling some result about potential theory and equilibrium measure and we enunciate our fundamental equilibrium measure which is the key of the work (theorem 2.4). In section 2 we prove the existence of such measure. Section 3 is dedicated to prove that the measure of theorem 2.4 is an equilibrium measure and we calculate explicitly the equilibrium energy of such measure.

   In section 4 we defined the ensemble of Generalized Gaussian random matrices ensemble and we prove the convergence of the statistical density of eigenvalues to the equilibrium measure of theorem 2.4. In the same way, we describe the probability
of atypical and large fluctuations of $\lambda_{\rm min}$ around its mean, say over a wider region of width $\sim O(\sqrt n)$. Since $\big<\lambda_{\rm min}\big>=b\sqrt{n}$, this requires the computation of the probability of an extremely rare event
characterizing a large deviation of $\sim O(\sqrt n)$ to the left of the mean. Such result has been proved in \cite{SD} in the case of Gaussian random matrices ensemble and our result follows in the same way.
\section{Solution on one single interval}
Let $\Sigma$ be a closed interval, and $Q$ be a lower semi-continuous function on $\Sigma$. If $\Sigma$ is unbounded we assume that
$$\lim_{|x|\to +\infty}\Big(Q(x)-\log(1+x^2)\Big)=\infty.$$
For given $Q$ and $\Sigma$, we wish to compute the equilibrium measure. We start by some general results.

For any probability measure $\mu$ on $\Sigma$, we defined the potential of $\mu$ by: for $x\in\Sigma\setminus{\rm supp}(\mu)$
$$U^\mu(x)=\int_\Sigma\log\frac{1}{|x-y|}\mu(dy),$$
and the energy by
$$E_{Q, \Sigma}(\mu)=\int_\Sigma U^\mu(x)\mu(dx)+\int_\Sigma Q(x)\mu(dx).$$
From the inequality $$|x-y|\leq\sqrt{1+x^2}\sqrt{1+y^2},$$
it can be seen that $E_{Q, \Sigma}(\mu)$ is bounded from below. Let
$$E^*_{Q, \Sigma}=\inf_{\mu\in\mathfrak{M}(\Sigma)}E_{Q, \Sigma}(\mu),$$
where $\mathfrak{M}(\Sigma)$ is the set of probability measures on the closed set $\Sigma$.

If $\mu(dx)=f(x)dx$, where $f$ is continuous function with compactly support $\subset\Sigma$, the potential is a continuous function, and $E_{Q, \Sigma}(\mu)<\infty$.
\begin{proposition}---
There is a unique measure $\mu^*\in\mathfrak{M}(\Sigma)$ such that $$E^*_{Q, \Sigma}=E_{Q, \Sigma}(\mu^*),$$
moreover the support of $\mu^*$ is compact.
\end{proposition}
This measure $\mu^*$ is called the {\it equilibrium measure.}

See theorem II.2.3 \cite{F}.
The next proposition is a method to find the equilibrium measure
\begin{proposition}---
Let $\mu\in\mathfrak{M}(\Sigma)$ with compact support.
Assume the potential $U^\mu$ of $\mu$ is continuous and, there is a
constant $C$ such that\\
{\rm (i)} $U^\mu(x)+\frac12Q(x)\geq C$ on $\Sigma$.\\
{\rm (ii)} $U^\mu(x)+\frac12Q(x)= C$ on {\rm supp}$(\mu)$. Then $\mu$ is the equilibrium measure $\mu=\mu^*$.

The constant C is called the (modified) Robin constant. Observe that $$E^*_{Q, \Sigma}=C+\frac12\int_\Sigma Q(x)\mu^*(dx).$$

\end{proposition}
\begin{remark}---
Let $\Sigma'\subset\Sigma$ be a closed interval of $\Sigma$, if we consider the restriction of the function $Q$ initially defined on $\Sigma$ to the closed interval $\Sigma'$ and if the equilibrium measure $\mu$ associate to $(\Sigma, Q)$ satisfies supp$(\mu)\subset\Sigma'$. Then the equilibrium measure for the couple $(\Sigma', Q)$ is $\mu$.
\end{remark}
We come to our first result. Let $\alpha\geq 0$, $\sigma>0$. Consider the closed interval $\Sigma_\sigma=[\sigma,+\infty[$ and
$$Q_\alpha(x)=x^2+2\alpha\log\frac1{x},$$
if $\alpha=0$, it is assumed that $\sigma\in\Bbb R$.

One can observe that $Q_\alpha$ is lower semi-continuous on the closed interval $\Sigma_\sigma$. Moreover $\lim\limits_{x\to +\infty}\Big(Q_\alpha(x)-\log(1+x^2)\Big)=+\infty$, hence the energy is correctly defined. {\it Let  $\nu_\alpha^\sigma$ be the equilibrium measure associate to $(\Sigma_\sigma, Q_\alpha)$}.
\newpage
\begin{theorem}---\
\begin{description}
 \item[(1)] If $\sigma=0$ and $\alpha>0$, then there is a unique $a_c=a(\alpha)>0$ and  a unique $b_c=b(\alpha)$ such that $b_c>\sqrt\alpha>a_c>0$, and the equilibrium measure on $\Sigma_0$ is given by the measure $\nu^0_\alpha$ with support $[a_c,b_c]$ and with density
$$f_{\alpha, 0}(x)=\frac1{\pi}\sqrt{(b_c-x)(x-a_c)}\Big(1+{\frac{\alpha}{\sqrt{a_cb_c}}}\frac1x\Big),$$
where $$\displaystyle b_c=\frac23\Big(\sqrt{6(\alpha+1)-2a_c^2}-\frac{a_c}2\Big).$$

\item[(2)] If $0<\sigma\leq a_c$ and $\alpha>0$, then the equilibrium measure on $\Sigma_\sigma$ still the same as in {\bf (1)} ($\nu^\sigma_\alpha=\nu^0_\alpha$).\\
 \item[(3)]If $\sigma>a_c$ and $\alpha>0$, in this case the equilibrium measure $\nu^\sigma_\alpha$ on $\Sigma_\sigma$ is supported by $]\sigma, b]$, and density $f_{\alpha, \sigma}$, where $b=b(\alpha, \sigma)$ is the unique solution of the following equations
$$\sigma+b-\frac{2\alpha}{\sqrt{\sigma b}}\geq 0,\;and\; \frac34(b-\sigma)^2+\sigma(b-\sigma)+2\alpha\sqrt{\frac \sigma b}-2\alpha-2=0,$$
and $$f_{\alpha, \sigma}(x)=\frac1{2\pi}\sqrt{\frac{b-x}{x-\sigma}}\Big(2x+b-\sigma-2\alpha\sqrt{\frac{\sigma}{b}}\frac1x\Big).$$

 \item[(4)]If $\alpha=0$ and $\sigma\in\Bbb R$, then two cases are present:
\begin{description}  \item[(a)] If $\sigma\geq-\sqrt 2$, the equilibrium measure $\nu^\sigma_0$ on $\Sigma_\sigma$ has support $]\sigma, b]$ and density $f_{0, \sigma}$, where $$b=b(0, \sigma)=\frac23\Big(\sqrt{\sigma^2+6}+\frac \sigma2\Big).$$
and $$f_{0, \sigma}(x)=\frac1{2\pi}\sqrt{\frac{b-x}{x-\sigma}}\Big(2x+b-\sigma\Big).$$
      \item[(b)] If $\sigma\leq-\sqrt 2$, the equilibrium measure $\nu_0^\sigma$ on $\Sigma_\sigma$ is the semicircle law with density $\displaystyle f_0(x)=\frac1{\pi}\sqrt{2-x^2}$, and support $[-\sqrt 2, \sqrt 2]$.
\end{description}
\end{description}
\end{theorem}

Before proving the theorem let gives same remarks.

 One can remark that, if $\sigma=0$ and $\alpha=0$, the equilibrium measure is supported by $]0,b_0]$ and have as density
$$f_{0, 0}(x)=\frac1{2\pi}\sqrt{\frac{b_0-x}{x}}\Big(2x+b_0\Big),$$
where $$b_0=b(0,0)=\frac23\sqrt6.$$
 Such density appear in the work (\cite{D}, \cite{S}) where the authors studies the density of positive eigenvalues of the Gaussian random matrices ensemble. It can be deduced from the first step. In fact $0<a_c<\sqrt\alpha$ and,  $$ b_c=b(a_c,\alpha)=\frac23\Big(\sqrt{6(\alpha+1)-2a_c^2}-\frac{a_c}2\Big),$$ hence as $\alpha\to 0$, we obtain $a_0=\lim\limits_{\alpha\to 0}a_c(\alpha)=0$ and $b_0=\lim\limits_{\alpha\to 0}b_c(\alpha)=\frac23\sqrt6.$\\
Moreover one can obtained such result from step {\bf (4)}.

The density in Step (4) appear in (\cite{D}, \cite{S}) where the authors studies the density of eigenvalues bigger then $\sigma$ in the case of the Gaussian random matrices ensemble. Such density can be obtained from step (3). Letting $\alpha\to 0$ one gets,\newline
$\sigma+b\geq 0$ and,  $$\frac34(b-\sigma)^2+\sigma(b-\sigma)-2=0.$$
The last equation has a unique solution $b$ with $b>\sigma$. Such solution can be find explicitly, it is given for all $\sigma\geq 0$ by $b=\frac23\Big(\sqrt{\sigma^2+6}+\frac \sigma2\Big).$

Moreover in this case, it will be seen that the measure $\nu^\sigma_0$ still an equilibrium measure for all $b\geq -\sigma$, and when $b=-\sigma$ which is the limit case, we obtain $b=\sqrt 2$, $\sigma=-\sqrt 2$, and $\nu^\sigma_0$ became the semi-circle law.\\
Observe that in {\bf(4)},{\bf (b)} the equilibrium measure is independent of the support $\Sigma_\sigma$. In fact for all $\sigma\leq-\sqrt 2$, the equilibrium measure relatively to the set $\Sigma_\sigma=[\sigma,+\infty[$, has semicircle law density.

To prove the theorem we need some preliminary results.
\section{Existence of the probability measure}
We will prove in the next proposition that, in each cases of the previous theorem, $\nu^\sigma_\alpha$ defined a probability measure. In section 3 we show that  such measure is an equilibrium measure.

For $\alpha>0$, let defined on $]0,+\infty[\times]0,+\infty[$,
$$\varphi_\alpha(x, a)=x+a-\frac{2\alpha}{\sqrt{ax}},$$
$$\psi_\alpha(x, a)=\frac{3}{4}(x-a)^2+a(x-a)+2\alpha{\frac{ \sqrt a}{\sqrt x}}-2\alpha-2,$$
and  $$E_\alpha=\Big\{a>0\mid\;\exists x>a,\;\varphi_\alpha(x, a)\geq 0\;{\rm and}\;\;\psi_\alpha(x, a)=0\Big\}.$$
For $\alpha=0$, $\varphi_\alpha$ and $\psi_\alpha$, will be defined on $\Bbb R^2$ and
$$E_0=\Big\{a\in\Bbb R\mid\;\exists x>a,\;\varphi_0(x, a)\geq 0\;{\rm and}\;\;\psi_0(x, a)=0\Big\}.$$

\begin{proposition}---
Let $\alpha\geq 0$, then
\begin{description}
\item [(1)]$E_\alpha$ is a closed set. Moreover $\displaystyle a_c=a_c(\alpha):=\min_aE_\alpha$ is correctly defined.
\item[(2)] If $\alpha>0$, then
\begin{description}
\item[(a)] for all $a> a_c$, there is a unique $b:=b(a ,\alpha)>a,$ such that $\psi_\alpha(b, a)=0$. Furthermore $a\mapsto b(a, \alpha)$ defined an increasing function and $\varphi_\alpha(b, a)\geq 0$
\item[(b)] For $a=a_c$ the unique element $b_c:=b(a_c, \alpha)>a_c$, satisfies $\varphi_\alpha(b_c, a_c)=0,\;\psi_\alpha(b_c, a_c)=0$ and it is given explicitly by
$$b_c=\frac23\Big(\sqrt{6(\alpha+1)-2a_c^2}-\frac{a_c}2\Big).$$
Moreover $0<a_c<\sqrt\alpha<b_c$.
\end{description}
\item[(3)] If $\alpha=0$, then $a_c=-\sqrt 2$, moreover,
\begin{description}
\item[(a)] if $a\geq-\sqrt2$, then the unique solution of the equation $\psi_0(x,a)=0$ in $]a,+\infty[$ is
$$b=\frac23\Big(\sqrt{a^2+6}+\frac a2\Big).$$
\item[(b)] If $a\leq-\sqrt2$, then $b=\sqrt 2$ is the unique solution of the equation $\psi_0(x, a)=0$ in $]a, +\infty[$.
\end{description}
\end{description}
\end{proposition}
Using mathematica we obtain, for $\alpha> 0$,
$$a_c=\sqrt{\frac53+\frac{5\alpha}{3}-\frac\gamma3-\frac23\sqrt{2+4\alpha-4\alpha^2+2(1+\alpha)\gamma}},$$
where $\displaystyle\gamma=\sqrt{1+2\alpha+4\alpha^2}.$

  To prove the proposition we need a technical lemma
 \begin{lemma}---Let $\beta$ be a real number, and $h$ be two times continuous differentiable function on $[\beta,+\infty[$ such that $h(\beta)<0$, $\lim\limits_{x\to +\infty}h(x)>0$, $\lim\limits_{x\to +\infty}h'(x)>0$ and $h''(x)>0$ for all $x\geq\beta$. Then there exist a unique $x_0>\beta$ such that $h(x_0)=0$.
\end{lemma}

 {\bf Proof.}--- The function $h'$ increase in the interval $[\beta,+\infty[$.\\
{\bf First case.} If $h'(\beta)\geq 0$, then $h'(x)\geq 0$, for all $x\geq\beta$. It follows that $h$ increase in $[\beta,+\infty[$. Since $h(\beta)<0$, and $\lim\limits_{x\to +\infty}h(x)>0$. Then the equation $h(x)=0$ admit a unique solution in the interval $]\beta,+\infty[$.\\
{\bf Second case.} If $h'(\beta)<0$. Since $\lim\limits_{x\to +\infty}h'(x)>0$, and $h'$ is an increasing function hence there exist a unique $\delta>\beta$, such that $h'(\delta)=0$, and by monotony of the function $h'$, $h'(x)<0$ in $[\beta,\delta[$, and, $h'(x)\geq 0$ in $[\delta,+\infty[$. Hence
$h$ decrease in $[\beta,\delta[$  and increase in $[\delta,+\infty[$. It follows that $h(x)\leq h(\beta)<0$ for all $x\in[\beta,\delta[$. Furthermore $\lim\limits_{x\to +\infty}h(x)>0$. Then there exist a unique $x_0\in[\delta,+\infty[$ such that $h(x_0)=0$. Which complete the proof.\\

{\bf Proof of proposition.}---\\
{\bf Step 1):} Let $\alpha>0$, and $a_n\in E_\alpha$ be a positive real sequence convergent to some $a$. Then by definition of $E_\alpha$ there is a sequence $x_n>a_n$, such that
\begin{equation}\label{12}\varphi_\alpha(x_n, a_n)\geq 0,\;{\rm and}\; \psi_\alpha(x_n, a_n)=0.\end{equation}
Since $\displaystyle\psi_\alpha(x_n, a_n)=\frac{3}{4}(x_n-a_n)^2+a_n(x_n-a_n)+2\alpha{\frac{ \sqrt a_n}{\sqrt x_n}}-2\alpha-2,$
hence $$a_n<x_n\leq \frac{2}{\sqrt  3}\sqrt{2\alpha+2}+a_n.$$
The sequence $a_n$ converge, it follows that $x_n$ is a bounded sequence, and there is some subsequence $x_{n_k}$ convergent to $x_0$, and $a\leq x_0$. From the inequality valuable for all $n$, $$\varphi_\alpha(x_n, a_n)=x_n+a_n-\frac{2\alpha}{\sqrt{a_nx_n}}\geq 0,$$ it follows that $\lim\limits_{n\to\infty}a_n=a>0$.
Using equation(\ref{12}), then when $n\rightarrow+\infty$, there is $x_0\geq a>0$, such that
$$\varphi_\alpha(x_0, a)\geq 0,\;{\rm and}\; \psi_\alpha(x_0, a)=0.$$
Since $\psi_\alpha(a, a)=-2$, hence $x_0>a$. which prove that $a\in E_\alpha$, and $E_\alpha$ is closed.

{\bf Existence of $a_c$: } Let $\alpha>0$. We will prove that $\sqrt\alpha\in E_\alpha$.
 Since $$\varphi_\alpha(\sqrt\alpha, \sqrt\alpha)=0.$$
Furthermore $x\mapsto\varphi_\alpha(x, \sqrt\alpha)$ is an increasing function hence $\varphi_\alpha(x, \sqrt\alpha)\geq 0$, for all $x\geq\sqrt\alpha$. So it is enough to prove that $\psi_\alpha(x, \sqrt\alpha)=0$ admit a solution $x_\alpha>\sqrt \alpha$. Take the derivative with respect to $x$, we obtain
$$\frac{\partial\psi_\alpha}{\partial x}(x, \sqrt\alpha)=\frac32(x-\sqrt\alpha)+\sqrt\alpha\,x-\frac{\alpha\sqrt{\sqrt\alpha}}{x\sqrt x}.$$
and $$\frac{\partial^2\psi_\alpha}{\partial^2 x}(x, \sqrt\alpha)=\frac32+\sqrt\alpha+\frac32\frac{\alpha\sqrt{\sqrt\alpha}}{x^2\sqrt x}.$$

Hence $\displaystyle\frac{\partial^2\psi_\alpha}{\partial^2 x}(x, \sqrt\alpha)>0$, for all $x\neq 0$. Furthermore $\psi_\alpha(\sqrt\alpha, \sqrt\alpha)=-2$, $\displaystyle\lim\limits_{x\to+\infty}\psi_\alpha(x, \sqrt\alpha)=\lim\limits_{x\to+\infty}\frac{\partial\psi_\alpha}{\partial x}(x, \sqrt\alpha)=+\infty$. Hence By the previous lemma, with $\beta=\sqrt\alpha$, it follows that for fixed $\alpha$, the equation $\psi_\alpha(x, \sqrt\alpha)=0$ admit a unique solution $x_\alpha$, in $]\sqrt\alpha, +\infty[$, which satisfies $\varphi_\alpha(x_\alpha, \sqrt\alpha)\geq 0$.

This prove that $E_\alpha\neq\varnothing$.

 For $\alpha>0$, $E_\alpha$ is bounded below by $0$, and by closeness $\min E_\alpha$ exist. It is obvious to see that $0\notin E_\alpha$ because $\varphi_\alpha(x, a)$ is not defined for $a=0$.\\
 For $\alpha=0$,
$$\varphi_0(x, a)=x+a,\;{\rm and}\;\psi_0(x, a)=\frac{3}{4}(x-a)^2+a(x-a)-2.$$
Hence for every $a\in\Bbb R$, the solutions of the equation $\psi_0(x, a)=0$  are
$$x_1=\frac{2}{3}\Big(\sqrt{a^2+6}+\frac a2\Big),\;{\rm and}\;x_2=-\frac{2}{3}\Big(\sqrt{a^2+6}-\frac a2\Big).$$
There is a unique solution $x\in]a, +\infty[$, such that $\varphi_0(x, a)\geq 0$. Such solution exist if $a\geq-\sqrt 2$, and is given by $$x=\frac{2}{3}\Big(\sqrt{a^2+6}+\frac a2\Big).$$
It follows that $a_c=\min_a E_\alpha=-\sqrt 2$.\\
{\bf Step(2):} Let $\alpha>0$. We begin by proving the existence, uniquness and the growth of  $a\mapsto b(a, \alpha)=x(a)$.\\
{\bf Case (a):} We saw that $a_c\in E_\alpha$, hence there exist some $b_c> a_c$, such that
\begin{equation}\label{op}
\varphi_\alpha(b_c, a_c)\geq 0,\;{\rm and}\;\psi_\alpha(b_c, a_c)=0.\end{equation}
For all $a>0$,
the function $x\mapsto \psi_\alpha(x, a)$ satisfies all the hypotheses of the previous lemma in the interval $[a,+\infty[$, in fact
$$\frac{\partial\psi_\alpha}{\partial x}(x, a)=\frac32(x-a)+ax-\frac{\alpha\sqrt a}{x\sqrt x},$$
and
$$\frac{\partial^2\psi_\alpha}{\partial^2 x}(x, a)=\frac32+a+\frac32\frac{\alpha\sqrt a}{x^2\sqrt x}>0,$$
$\psi_\alpha(a, a)=-2$, $\displaystyle\lim\limits_{x\to +\infty}\psi_\alpha(x, a)=\lim\limits_{x\to +\infty}\frac{\partial\psi_\alpha}{\partial x}(x,a)=+\infty$.
it follows that,
 \begin{equation}\label{oq}\exists\; !x\in ]a,+\infty[;\;\;\psi_\alpha\big(x(a), a\big)=0.\end{equation}
By equations (\ref{op}), (\ref{oq}) the set $\Big\{x(a)\mid a>0\Big\}$ where $x(a)$ is a the unique solution in $]a,+\infty[$ of the problem
$$\varphi_\alpha(x, a)\geq 0,\;{\rm and}\;\psi_\alpha(x, a)=0$$ is not empty. By unicity, it follows that $a\mapsto x(a)$ defined a function on $]0,+\infty[$.

{\bf Growth of the function $a\mapsto x(a)$.}
Let $0<a_1<a_2$, and $x(a_1)$ be a solution of the problem
$$\exists\; x>a_1>0,\;\varphi_\alpha\big(x, a_1\big)\geq 0,\;\psi_\alpha\big(x, a_1\big)=0,$$
and $x(a_2)$ the unique element in the interval $]a_2, +\infty[$, such that
$$\psi_\alpha\big(x(a_2), a_2\big)=0.$$
We want to prove $x(a_1)\leq x(a_2)$. Assume the contraire. Then
$$a_1<a_2< x(a_2)<x(a_1).$$
Take the derivative of the function $a\mapsto \psi_\alpha(x, a)$, it yields
$$\frac{\partial\psi_\alpha}{\partial a}(x,a)=-\frac12\varphi_\alpha(x,a),$$
and $$\frac{\partial^2\psi_\alpha}{\partial^2 a}(x,a)=-\frac12(1+\frac{\alpha}{a\sqrt{ax}})<0.$$
Hence for all $x>0$, the function $\displaystyle a\mapsto\frac{\partial\psi_\alpha}{\partial a}(x,a)$ decreases on $]0,+\infty[$. Thus for all $a>a_1$,
$$\frac{\partial\psi_\alpha}{\partial a}\big(x(a_1),a\big)\leq\frac{\partial\psi_\alpha}{\partial a}\big(x(a_1),a_1\big)=-\frac12\varphi_\alpha\big(x(a_1),a_1\big)\leq 0,$$
which mean that the function $a\mapsto \psi_\alpha\big(x(a_1),a\big)$ decreases on $]a_1,+\infty[$.

Since $\displaystyle a_1<a_2,$ it follows that
$$\psi_\alpha\big(x(a_1),a_2\big)\leq \psi_\alpha\big(x(a_1),a_1\big)=0.$$
By the assumption $a_2< x(a_2)<x(a_1),$ and the unicity of the solution of the equation $\psi_\alpha(x,a_2)=0$ in $]a_2,+\infty[$, on gets
\begin{equation}\label{oo}\psi_\alpha\big(x(a_1),a_2\big)<0.\end{equation}
Now, consider the function $x\mapsto \psi_\alpha(x,a_2)$. Such function satisfies for all $x>0$,
$$\frac{\partial^2\psi_\alpha}{\partial^2 x}(x,a_2)>0,$$
then $\displaystyle x\mapsto \frac{\partial\psi_\alpha}{\partial x}(x,a_2)$ increases in $]a_2,+\infty[$.\\
If  $\displaystyle \frac{\partial\psi_\alpha}{\partial x}(a_2,a_2)\geq 0$. Then  $\displaystyle \frac{\partial\psi_\alpha}{\partial x}(x,a_2)\geq 0$, for all $x\geq a_2$, and the function $\displaystyle x\mapsto \psi_\alpha(x,a_2)$ increases in $]a_2,+\infty[$. Using the assumption  $x(a_2)<x(a_1)$ and equation(\ref{oo}), one gets
$$0=\psi_\alpha\big(x(a_2), a_2\big)\leq\psi_\alpha(x(a_1),a_2)<0.$$
this gives a contradiction.\\
Assume $\displaystyle\frac{\partial\psi_\alpha}{\partial x}(a_2,a_2)<0$. The function $\displaystyle x\mapsto \frac{\partial\psi_\alpha}{\partial x}(x,a_2)$ is strictly increasing on $]a_2, +\infty[$, and $\displaystyle \lim\limits_{x\to +\infty} \frac{\partial\psi_\alpha}{\partial x}(x,a_2)=+\infty$, hence there is a unique $\delta>a_2$, such that $\displaystyle\frac{\partial\psi_\alpha}{\partial x}(\delta,a_2)=0.$ It follows that $x\mapsto \psi_\alpha(x, a_2)$ decreases on $[a_2, \delta]$, and increases on $[\delta, +\infty[$. Moreover $\displaystyle \psi_\alpha(a_2, a_2)=-2$, hence $\displaystyle x(a_2)\in]\delta, +\infty[$. From the assumption
$x(a_2)<x(a_1)$, and equation(\ref{oo}) we obtains
$$0=\psi_\alpha(x(a_2), a_2)\leq\psi_\alpha(x(a_1), a_2)<0.$$
This gives contradiction. Which is the desired result.

Let $\sigma> a_c$. From the previous, there is a unique $b=b(\sigma, \alpha)>\sigma$, such that $\psi_\alpha(b, \sigma)=0$. Moreover $\varphi_\alpha(b_c,a_c)\geq 0$ and $\psi_\alpha(b_c,a_c)=0$, and the function $\sigma\mapsto b(\sigma, \alpha)$ increase, hence $b=b(\sigma,\alpha)\geq b_c=b(a_c,\alpha)$. Furthermore the two functions $\displaystyle x\mapsto\varphi_\alpha(x,a)$, and $\displaystyle a\mapsto\varphi_\alpha(x,a)$ increases. Thus
$$\varphi_\alpha(b,\sigma)\geq\varphi_\alpha(b_c,\sigma)\geq\varphi_\alpha(b_c,a_c)\geq 0.$$ Which means that $b$ is the unique solution in $]\sigma,+\infty[$ of the problem
$$\varphi_\alpha(x,\sigma)\geq 0,\;\psi_\alpha(x,\sigma)=0.$$

{\bf (b)}  Assume $\sigma=a_c$. We saw that $a_c>0$, hence there exist $n_0$, such that for all $n\geq n_0$, $a_n=a_c-\frac1n>0$. By the definition of $a_c$, for all $x\geq a_n$, $\varphi_\alpha(x, a_n)<0$, or $\psi_\alpha(x, a_n)\neq 0$. Furthermore the equation $\psi(x, a)=0$ has a unique solution in $]a,+\infty[$, for each $a>0$. It follows that for all $n\geq n_0$, there is some sequence $x_n$ such that $x_n>a_n$, $\psi_{\alpha}(x_n, a_n)= 0$ and $\varphi_{\alpha}(x_n, a_n)<0$.

 Since $$\varphi_{\alpha}(x_n, a_n)=x_n+a_n-\frac{2\alpha}{\sqrt{a_nx_n}}<0,$$
 Hence $$a_n<x_n\leq -a_n+\frac{2\alpha}{a_n}.$$
 The sequence $a_n$ converge to $a_c$, it follows that $x_n$ is bounded, and there is some subsequence $x_{n_k}$ convergent to $x_0$, and $\displaystyle 0<a_c\leq x_0\leq -a_c+\frac{2\alpha}{a_c}.$
 Since $$\psi_{\alpha}(x_{n_k}, a_{n_k})=0,\;\;{\rm and}\;\varphi_{\alpha}(x_{n_k}, a_{n_k})<0.$$
 letting $k$ to infinity, we obtain by continuity
 $$\psi_{\alpha}(x_0, a_c)=0,\;\;{\rm and}\;\varphi_{\alpha}(x_0, a_c)\leq 0.$$
 By unicity of the solution of the equation $\psi_{\alpha}(x, a_c)=0$ in the interval $]a_c,+\infty[$, it follows that $x_0=b_c$. Hence $\varphi_{\alpha}(b_c, a_c)\leq 0,$ and by the definition of $b_c$, it follows that $\varphi_{\alpha}(b_c, a_c)= 0$ and $\psi_\alpha(b_c, a_c)=0$. Which means that $$b_c+a_c-\frac{2\alpha}{\sqrt{a_cb_c}}=0,$$
 and $$\frac34(b_c-a_c)^2+a_c(b_c-a_c)+2\alpha\frac{\sqrt{a_c}}{\sqrt{b_c}}-2\alpha-2=0,$$
 let put $y=b_c-a_c$, it yield from the two previous equations
  $$\frac34y^2+2a_cy+2a_c^2-2\alpha-2=0.$$
  By simple computation we obtain $$y=\frac43\Bigg(-a_c+\sqrt{\frac32(\alpha+1)-\frac{a_c^2}{2}}\Bigg).$$
 Observe that $y>0$ and the square root is correctly defined. In fact
  $a_c<b_c$, hence $$2a_c<a_c+b_c=\frac{2\alpha}{\sqrt{a_cb_c}}<\frac{2\alpha}{a_c},$$
  and $$\frac{2\alpha}{b_c}<\frac{2\alpha}{\sqrt{a_cb_c}}=a_c+b_c<2b_c$$
  it yield $$0<a_c<\sqrt\alpha<b_c,$$
  and The value of $b_c$ is given by $$b_c=\frac23\Big(\sqrt{6(\alpha+1)-2a_c^2}-\frac{a_c}{2}\Big).$$
{\bf Step 3):} Assume $\alpha=0$, then\\
  $\displaystyle\psi_\alpha(x, a)=\frac34(x-a)^2+a(x-a)-2=0$ and $\displaystyle\varphi_\alpha(x, a)=x+a$. The solutions of the equation  $\psi_\alpha(x, a)=0$ are known explicitly
 $$x_1(a)=\frac23\big(-a-\sqrt{a^2+6}\big)+a,\quad x_2(a)=\frac23\big(-a+\sqrt{a^2+6}\big)+a,$$
 but the two inequalities $\varphi_\alpha(x, a)\geq 0$, and $x>a$ hold if and only if
 $\displaystyle x=\frac23\big(-a+\sqrt{a^2+6}\big)+a$ and $a\geq -\sqrt 2$.
  Then the unique solution in the interval $]a,+\infty[$ with $a\geq-\sqrt 2$ is
 $$b=b(a)=\frac23(\sqrt{a^2+6}+\frac a2).$$
This complete the proof.\\

 \begin{center}{\bf Graphics}\end{center} Here we plot the two functions $x\mapsto\varphi_\alpha(x,a)$ and $x\mapsto\psi_\alpha(x,a)$, for the same value $\alpha=2$, and different values of $a$.\\

\newpage
\begin{figure}[h]
\centering\scalebox{0.5}{\includegraphics[width=14cm, height=12cm]{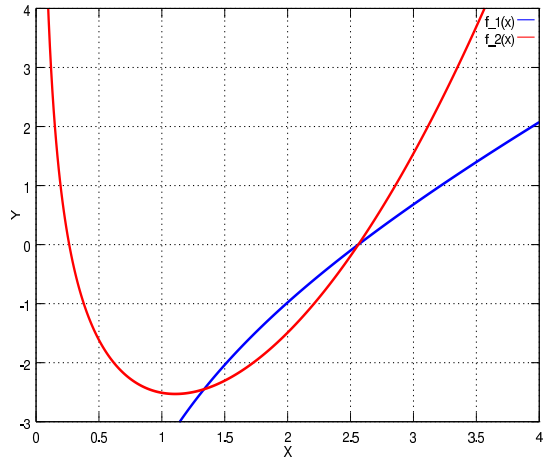}}
\caption{\bf $\alpha=2$, $a\approx a_c$.}
\end{figure}
\begin{flushleft}\textcolor{blue}{\rule{1cm}{1pt}} $f_1(x)=\varphi_\alpha(x,a_c)$.\\
\textcolor{red}{\rule{1cm}{1pt}} $f_2(x)=\psi_\alpha(x,a_c)$.\\ $a_c\approx 0.618$, $b_c\approx 2.562$, $\varphi_\alpha(b_c,a_c)=0$, $\psi_\alpha(b_c,a_c)=0$.
\end{flushleft}

\begin{figure}[h]
\centering\scalebox{0.5}{\includegraphics[width=14cm, height=12cm]{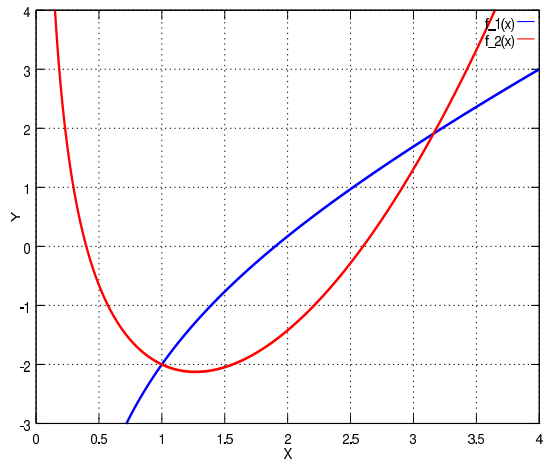}}
\caption{\bf $\alpha=2$, $a_c<a$.}
\end{figure}
\begin{flushleft}
\textcolor{blue}{\rule{1cm}{1pt}} $f_1(x)=\varphi_\alpha(x,a)$.\\
\textcolor{red}{\rule{1cm}{1pt}} $f_2(x)=\psi_\alpha(x,a)$.\\
$a=1>a_c$, $b=b(a, \alpha)\approx 2.6$ and $\varphi_\alpha(b, a)=1.11>0$, $\psi_\alpha(b, a)=0$.
\end{flushleft}
\newpage
\begin{figure}[h]
\centering\scalebox{0.5}{\includegraphics[width=14cm, height=12cm]{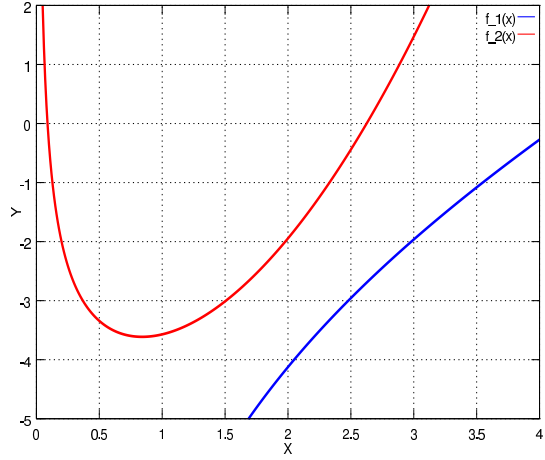}}
\caption{\bf $\alpha=2$, $a<a_c$.}

\end{figure}
\begin{flushleft}
\textcolor{blue}{\rule{1cm}{1pt}} $f_1(x)=\varphi_\alpha(x, a)$.\\
\textcolor{red}{\rule{1cm}{1pt}} $f_2(x)=\psi_\alpha(x, a)$.\\ $a=0.2<a_c$, $b=b(a, \alpha)\approx 2.62 $, $\varphi_\alpha(b, a)\approx -2.67<0$ and $\psi_\alpha(b,a)=0$.
\end{flushleft}
\section{Proof of theorem1.4}
For $\alpha\geq 0$, $a>0$ and $b>a$, let $\nu_\alpha^a$ be the measure supported by $]a,b]$ and with density
$$f_{\alpha, a}(x)=\frac{1}{2\pi}\sqrt{\frac{b-x}{x-a}}\Big(2x+b-a-\alpha\sqrt{\frac{a}{b}}\frac1x\Big).$$
If $\alpha=0$, it will be assumed that $a\in\Bbb R$.
In this section we will prove that $\nu_\alpha^\sigma$ is an equilibrium measure.
\begin{lemma}---
For $\alpha\geq 0$ the measure $\nu^a_\alpha$, with density $f_{\alpha, a}$ defined a probability measure on $]a,b]$ if and only if $(a, b)$ is the unique solution given in proposition 2.1.
\end{lemma}
 {\bf Proof.}---Let $b>a$. Then $\nu^a_\alpha$ is a probability if and only if the function $f_{\alpha, a}$ is positive with integral one over $]a,b]$.
Since
$$\int_a^bf_{\alpha,a }(x)dx=\frac1{\pi}\int_a^bx\sqrt{\frac{b-x}{x-a}}dx+\frac{b-a}{2\pi}\int_a^b\sqrt{\frac{b-x}{x-a}}
dx-\frac{\alpha}{\pi}\sqrt{\frac{a}{b}}
\int_a^b\sqrt{\frac{b-x}{x-a}}\frac{dx}{x},$$
performing the change of variable $u=\frac{x-a}{b-a}$ for the first and the second integral, it follows that

$$I=\int_a^b\sqrt{\frac{b-x}{x-a}}dx=(b-a)\int_0^1\sqrt{\frac{1-x}{x}}=\pi\frac{b-a}{2},$$
and
$$\int_a^bx\sqrt{\frac{b-x}{x-a}}dx=a I+(b-a)^2\int_0^1\sqrt{(1-u)u}du=\pi\frac{a(b-a)}{2}
+\frac\pi8(b-a)^2.$$
The last integral can be written as
$$\int_a^b\sqrt{\frac{b-x}{x-a}}\frac{dx}{x}=\frac1a\Big(I-\int_a^b\sqrt{(b-x)(x-a)}
\frac{dx}{x}\Big).$$
Since $\displaystyle\int_a^b\sqrt{(b-x)(x-a)}
\frac{dx}{x}=\frac{\pi}{2}(\sqrt b-\sqrt a)^2$ (see the Appendix), hence
$$\int_a^b\sqrt{\frac{b-x}{x-a}}\frac{dx}{x}=\frac{\pi}{\sqrt a}(\sqrt{b}-\sqrt a).$$
Combining all these results, togethers gives
$$\int_a^bf_{\alpha, a}(x)dx=\frac18(b-a)^2+\frac{(b-a)^2}{4}+\frac{a(b-a)}{2}
-\alpha\sqrt{\frac{a}{b}}(\sqrt{\frac{b}{a}}-1).$$
Hence $\nu^a_\alpha$ is a probability measure, if and only if,
$$\frac34(b-a)^2+a(b-a)+2\alpha\frac{{\sqrt a}}{{\sqrt b}}-2\alpha-2=0\quad{\rm and}\; f_{\alpha, a}(x)\geq 0,\;\forall\;x\in]a,b].$$
Since $f_{\alpha, a}(x)\geq 0$ if and only if, for all $x\in]a,b]$, $$\displaystyle h(x)=2x+b-a-2\alpha\sqrt{\frac{a}{b}}\frac1x\geq 0.$$
The function $h$ is continuous and increasing. Hence $h(x)\geq 0$, for all $x\in]a,b]$, if and only if $h(a)\geq 0$, which means that
$$b+a-\frac{2\alpha}{\sqrt{a b}}\geq 0.$$
Hence $\nu^a_\alpha$ is a probability measure if and only if $a\geq a_c$, and  $b>a$ is the unique solution of
$$\psi_\alpha(b, a)=0\;\;{\rm and}\;\;\varphi_\alpha(b, a)\geq 0.$$
From proposition 2.1, such solution exist and is unique.

 This complete the proof.

\subsection{Logarithm potential of $\nu^a_\alpha$}
For $\alpha\geq 0$, let $(a,b)$ be a pair defined as in proposition 2.1.
The Cauchy transform $G_{\nu^a_\alpha}$ of the measure $\nu^a_\alpha$, defined on $\Bbb C\setminus[a,b]$ by
$$G_{\nu^a_\alpha}(z)=\int_a^b\frac{1}{z-x}\nu^a_\alpha(dx),$$
The function $\displaystyle g_{\alpha, a}(\omega)=\sqrt{\frac{\omega-b}{\omega-a}}\Big(\omega+\frac{b-a}{2}-\alpha\sqrt{\frac{a}{b}}\frac1\omega\Big),$ is meromorphic on $\Bbb C\setminus[a,b]$ and the boundary values of  $g_{\alpha, a}$ is given by
$$[g_{\alpha, a}]=2i\pi f_{\alpha, a}(x)\chi(x)+2i\pi\alpha\delta_0.$$
where $\chi$ is the indicator function of the interval $[a,b]$.

Furthermore as $|\omega|\to+\infty$, $$g_{\alpha, a}(\omega)=\omega+\frac{c}{\omega}+O(\frac{1}{\omega}),$$

the residues of $\displaystyle \frac{1}{z-\omega}g_{\alpha, a}(\omega)$ are: on $\omega=z$ is given by $-g_{\alpha, a}(z)$, on $0$ is $-\displaystyle \frac{\alpha}{z}$ and the residues at infinity is $z$.\\
It follows from the Residues and Liouville  theorems that
$$G_{\nu^a_\alpha}(z)=-g_{\alpha, a}(z)+z-\frac{\alpha}{z}.$$

For $\sigma\in\Bbb R$, let $Q$ be a continuous function on $\Sigma_\sigma=[\sigma,+\infty[$. We defined the potential of a probability measure $\mu$ with compact support in $\Sigma_\sigma$, by
 $$U^{\mu}(x)=\int_\sigma^{+\infty}\log\frac{1}{|x-t|}\mu(dt),$$
 and the energy $E_{Q, \sigma}(\mu)$ by
 $$E_{Q, \sigma}(\mu)=\int_\sigma^{+\infty}U^\mu(x)\mu(dx)+\int_\sigma^{+\infty}Q(x)\mu(dx).$$
 \begin{proposition}---For all $\sigma\geq 0$, and $\alpha\geq 0$, let consider on the interval $\Sigma_\sigma$, the function
 $\displaystyle Q_{\alpha}(x)=x^2+2\alpha\log\frac1x$. Then there exist some constant $C$, such that
 \begin{description}
 \item[(i)] $U^{\nu^a_\alpha}(x)+\frac12Q_{\alpha}(x)=C$\;\;on $]a,b]$.
 \item[(ii)] $U^{\nu^a_\alpha}(x)+\frac12Q_{\alpha}(x)\geq C$\;\;on $[\sigma,a]\cup[b,+\infty[$.

 \end{description}
 \end{proposition}
  It follows from the proposition that, the equilibrium measure which minimize the energy is $\nu^a_\alpha$, and the equilibrium energy is given by
 $$\displaystyle E^*_{\alpha, \sigma}=C+\frac12\int_a^bQ_{\alpha}(x)\nu^{a}_\alpha(dx).$$
 {\bf Proof}---.We saw that
 $$G_{\nu^a_\alpha}(z)=z-\frac{\alpha}{z}-g_{\alpha, a}(z).$$
 {\bf (1):} $\sigma=0$, $\alpha>0$, the minimum value of $a$ such that $\nu^a_\alpha$ is positive is $a=a_c$ and to obtain a probability measure, $b$ must be $b=b_c$
 It follows that:\\
 for $0<x\leq a_c$, $${\rm Re}G_{\nu^a_\alpha}(x)=x-\frac{\alpha}{ x}+h_{\alpha, a}(x),$$
 for $a_c\leq x\leq b_c$, $${\rm Re}G_{\nu^a_\alpha}(x)=x-\frac{\alpha}{ x},$$
 for $x\ge b_c$,$${\rm Re}G_{\nu^a_\alpha}(x)=x-\frac{\alpha}{ x}-h_{\alpha, a}(x).$$
where $$h_{\alpha, a}(x)=\sqrt{(b_c-x)(a_c-x)}\left(1+\frac{\alpha}{\sqrt{a_cb_c}}\frac1x\right).$$
 Since $$\frac{d}{dx}U^\mu(x)=-{\rm Re}G_{\nu^a_\alpha}(x),$$
 it follows that
  $$\begin{aligned}
 U^{\nu^a_\alpha}(x)&=-\Big(\frac12x^2+\alpha\log\frac1x\Big)+C+\int_x^{a_c}h_{\alpha, a}(t)\,dt\quad {\rm if}\;0< x\leq a_c,\\
 U^{\nu^a_\alpha}(x)&=-\Big(\frac12x^2+\alpha\log\frac1x\Big)+ C\quad\quad\quad\quad\quad\quad\quad{\rm if}\;a_c\leq x\leq b_c,\\
 U^{\nu^a_\alpha}(x)&=-\Big(\frac12x^2+\alpha\log\frac1x\Big)+C+\int_{b_c}^xg_{\alpha, a}(t)\,dt\;\;\quad {\rm if}\;x\geq b_c.\\
 \end{aligned}$$
  Hence
 $$ U^{\nu^a_\alpha}(x)+\frac12Q_\alpha(x)\left\{\begin{aligned}&=C\quad {\rm on}\;[a_c,b_c]\\&\geq C\;\;{\rm on}\;]\sigma,a_c]\cup[b_c,+\infty[\end{aligned}\right.$$
{\bf (2):} $0<\sigma\leq a_c$, $\alpha>0$, here also the unique values such that $\nu^a_\alpha$ is a probability are $a=a_c$, and $b=b_c$, and the result is as in {\bf (1)}. This result can be derived from remark 1.3\\
{\bf (3)} $a_c<\sigma$, $\alpha>0$, in this case for $\nu^a_\alpha$ to be a probability with support $\subset [\sigma,+\infty[$, we must have $a=\sigma$, for this value of $a$, there is a unique $b=b(\sigma,\alpha)$ as in proposition 2.1. Moreover,\\
for $\sigma\leq x\leq b$, $${\rm Re}G_{\nu^a_\alpha}(x)=x-\frac{\alpha}{ x},$$
 for $x\geq b$,$${\rm Re}G_{\nu^a_\alpha}(x)=x-\frac{\alpha}{ x}-g_{\alpha, a}(x).$$
 where $$g_{\alpha, \sigma}(x)=\sqrt{\frac{x-b}{x-\sigma}}\left(x+\frac{b-\sigma}{2}-\alpha\sqrt{\frac{\sigma}{b}}\frac1x\right),$$
 $$\begin{aligned}
 U^{\nu^a_\alpha}(x)&=-\Big(\frac12x^2+\alpha\log\frac1x\Big)+ C\quad\quad\quad\quad\quad\quad\;\;\;\;{\rm if}\;\sigma\leq x\leq b,\\
 U^{\nu^a_\alpha}(x)&=-\Big(\frac12x^2+\alpha\log\frac1x\Big)+C+\int_{b}^xg_{\alpha, a}(t)\,dt\quad {\rm if}\;\;\;x\geq b.\\
 \end{aligned}$$
 Since from proposition 2.1, for all $x>b$, $$g_{\alpha, \sigma}(x)\geq\sqrt{\frac{x-b}{x-\sigma}}\varphi_\alpha(b, \sigma)\geq 0,$$ hence
 $$ U^{\nu^a_\alpha}(x)+\frac12Q_\alpha(x)\left\{\begin{aligned}&=C\quad {\rm on}\;[\sigma,b]\\&\geq C\;\;{\rm on}\; [b,+\infty[\end{aligned}\right.$$
 {\bf (4):} $\alpha=0$ and $\sigma\in\Bbb R$, here if $\sigma\geq-\sqrt 2$ we put $a=\sigma$, and if $\sigma<-\sqrt 2$, we put $a=-\sqrt 2$, hence as in proposition 2.1, $\displaystyle b(a)=\frac23\Big(\sqrt{a^2+6}+\frac{a}{2}\Big)$, and $b(-\sqrt 2)=\sqrt2$. In the two cases $\nu^a_0$ is a probability.\\
 {\bf case 1:}  $\sigma\geq-\sqrt 2$, hence $a=\sigma$ and $b=\frac23\Big(\sqrt{a^2+6}+\frac{a}{2}\Big)$. Observe that $\varphi_\alpha(b, a)=b+a\geq 0$, moreover:\\
 for $a\leq x\leq b$, $${\rm Re}G_{\nu^a_0}(x)=x,$$
 for $x\geq b$,$${\rm Re}G_{\nu_0}(x)=x-g_{0, a}(x),$$
 where $$g_{0, a}(x)=\sqrt{\frac{x-b}{x-a}}\Big(x+\frac{b-a}{2}\Big),$$
 hence
 $$\begin{aligned}
 U^{\nu^a_0}(x)&=-\frac12x^2+ C\quad\quad\quad\quad\quad\quad\;\;{\rm if}\;a\leq x\leq b,\\
 U^{\nu^a_0}(x)&=-\frac12x^2+C+\int_{b}^xg_{0, a}(t)\,dt\; \quad{\rm if}\;x\geq b.\\
 \end{aligned}$$
 and
 $$ U^{\nu^a_0}(x)+\frac12Q_0(x)\left\{\begin{aligned}&=C\quad {\rm on}\;[a,b]\\&\geq C\;\;{\rm on}\; [b,+\infty[\end{aligned}\right.$$
 {\bf case 2:} $\sigma<-\sqrt 2$, hence $a=-\sqrt 2$, and $b=\sqrt 2$:\\
 for $\sigma\leq x\leq -\sqrt 2$,$${\rm Re}G_{\nu^a_0}(x)=x+\sqrt{x^2-2},$$
 for $-\sqrt 2\leq x\leq \sqrt 2$, $${\rm Re}G_{\nu^a_0}(x)=x,$$
 for $x\geq \sqrt 2$,$${\rm Re}G_{\nu^a_0}(x)=x-\sqrt{x^2-2}.$$
 Hence $$\begin{aligned}
 U^{\nu^a_0}(x)&=-\frac12x^2+C+\int_x^{-\sqrt 2}\sqrt{t^2-2}\,dt\quad {\rm if}\;\sigma\leq x\leq -\sqrt 2,\\
 U^{\nu^a_0}(x)&=-\frac12x^2+ C\quad\quad\quad\quad\quad\quad\quad\quad\;\;\;{\rm if}\;-\sqrt 2\leq x\leq\sqrt 2,\\
 U^{\nu^a_0}(x)&=-\frac12x^2+C+\int_{\sqrt 2}^x\sqrt{t^2-2}\,dt\quad\quad {\rm if}\;x\geq\sqrt 2.\\
 \end{aligned}$$
 and
 $$ U^{\nu^a_0}(x)+\frac12Q_0(x)\left\{\begin{aligned}&=C\quad {\rm on}\;[-\sqrt 2,\sqrt 2]\\&\geq C\;\;{\rm on}\; [\sigma, -\sqrt 2]\cup[\sqrt 2,+\infty[\end{aligned}\right.$$
 which complete the proof.
 \subsection{Computation of the energy $E^*_{\alpha, \sigma}$}
 We saw that $$E^*_{\alpha, \sigma}=C+\frac12\int_a^bQ_\alpha(x)\nu^a_\alpha(dx),$$
 we begin by computing the modified Robin constant.
 \begin{proposition}---
 The Modified Robin constant is given by
 $$\begin{aligned}C&=-\frac{1}{16} \Big(a^2+6 a b+b^2+2 (a-b)^2 \log\frac{b-a}{4}-4(a+b)^2+4(b^2-a^2)\log\frac{b-a}{4}\Big)\\&- {\frac{\alpha}{\sqrt {ab}}}  \left(\frac{a+b}{2}-\frac{b-a}{2}\log\frac{b-a}{4}-(a+b)(\frac12+\log 2-\frac12\log(b-a))+\sqrt{ab}\log{\frac{\sqrt b+\sqrt a}{\sqrt b-\sqrt a}}\right)\end{aligned}$$
  where $b>a>0 $, which are given as in proposition 2.1.
 \end{proposition}
  One observe that when $\alpha=0$ and $\sigma=0$, hence $a=0$ and $b^2=\frac83$ moreover $$C=\frac12+\frac12\log2+\frac12\log3.$$ We recover's the (modified) Robin constant of Dean-Majumdar distribution.
 \begin{lemma}---
 For all $0<a<b$, Let $\mu$ be the positive measure defined on $[0,+\infty[$ by
 $$\int_0^{+\infty}f(t)\mu(dt)=\frac1{2\pi}\int_a^bf(t)\sqrt{(t-a)(b-t)}\frac{dt}{t}.$$
 Then
 $$\frac1{2\pi}\int_a^b\sqrt{(t-a)(b-t)}\frac{dt}{t}=\frac14\big(\sqrt b-\sqrt a\big)^2.$$
 And
 $$\lim_{x\to+\infty}U^\mu(x)+\frac14\Big(\sqrt b-\sqrt a\Big)^2\log x=0.$$
 Moreover
 the logarithmic potential $U^\mu$ is given for all $x\geq b$, by
 $$U^\mu(x)=-\frac12(x-\sqrt {ab}\log x)+C_\mu+\frac12\int_b^x\sqrt{(t-a)(t-b)}\frac{dt}{t},$$
 where $\displaystyle C_\mu=\frac14(a+b)\Big(1-\log\frac{b-a}{4}\Big)-\frac{\sqrt{ab}}{2}\log{\frac{\sqrt b+\sqrt a}{\sqrt b-\sqrt a}}$.
 \end{lemma}
 The lemma will be proved in the Appendix.

{\bf Proof of proposition 3.3.}---To compute the constant $C$, we will use the fact that $\displaystyle \lim_{x\to+\infty}U^{\nu^a_\alpha}(x)+\log x=0$. Since for $x\geq b$,
\begin{equation}\label{011} U^{\nu^a_\alpha}(x)=-(\frac12x^2+\alpha\log\frac1x)+C+\int_b^xg_{\alpha, a}(t)\,dt\end{equation}
 Let now compute the integral of $g_{\alpha, a}$.
 $$\int_b^xg_{\alpha, a}(t)\,dt=\int_b^x \sqrt{\frac{t-b}{t-a}}\Big(t+\frac{b-a}{2}-\alpha\sqrt{\frac{a}{b}}\frac1t\Big)dt,$$
 hence,
$$\int_b^xg_{\alpha, a}(t)\,dt= \int_b^x \sqrt{(t-b)(t-a)}dt+\frac{a+b}{2}\int_b^x \sqrt{\frac{t-b}{t-a}}dt-
\alpha\sqrt{\frac{a}{b}}\int_b^x\sqrt{\frac{t-b}{t-a}}\frac{dt}{t}.$$
Since $\displaystyle \sqrt{\frac{t-b}{t-a}}\frac1t=\frac1a\Big(\sqrt{\frac{t-b}{t-a}}-\frac1t\sqrt{(t-b)(t-a)}\Big),$ it follows that
$$\begin{aligned}\int_b^xg_{\alpha, a}(t)\,dt&= \int_b^x \sqrt{(t-b)(t-a)}dt+\Big(\frac{a+b}{2}-\frac{\alpha}{\sqrt{ab}}\Big)\int_b^x \sqrt{\frac{t-b}{t-a}}dt\\&+
\frac{\alpha}{\sqrt{ab}}\int_b^x\sqrt{(t-b)(t-a)}\frac{dt}{t}.\end{aligned}$$
 Let
 $$I_1(x)=\int_b^x \sqrt{(t-b)(t-a)}dt,$$
 $$I_2(x)=\int_b^x \sqrt{\frac{t-b}{t-a}}dt,$$
 $$I_3(x)=\int_b^x\sqrt{(t-b)(t-a)}\frac{dt}{t}.$$
Since
 \begin{equation}\label{00}I_1(x)=\frac12x^2-\frac{a+b}{2}x-\frac{\theta^2}{8}\log x+A_1+o(1),\end{equation}
 where $$\displaystyle A_1=\frac1{16}\Big(a^2+6ab+b^2+2(b-a)^2\log\frac{b-a}{4}\Big).$$
 Moreover
 \begin{equation}\label{000}I_2(x)=x-\frac\theta 2\log x+A_2+o(1),\end{equation}
 where $\displaystyle A_2= -\frac{a+b}{2}+\frac{b-a}{2} \log\frac{b-a}{4}.$
  See lemma6.2 in the appendix.\\
  Using the value of $I_3$ in the previous lemma, and equations (\ref{00}),(\ref{000}) we obtain
 $$\int_b^xg_{\alpha, a}(t)\,dt=\frac12x^2+\Big(-\frac{\theta^2}{8}-\frac{\alpha(a+b)}{2\sqrt{ab}}+\frac{\alpha\theta}{2\sqrt{ab}}-\frac{\theta(a+b)}{4}\Big)\log x+A+o(\frac1x).$$
 where $\displaystyle A=A_1+A_2(\frac{a+b}{2}-\frac{\alpha}{\sqrt{ab}})+\frac{\alpha}{\sqrt{ab}}C_\mu.$
 Furthermore, by making use of proposition 2.1 $(\psi_\alpha(b, a)=0)$, one gets
 $$\begin{aligned}&-\frac{\theta^2}{8}-\frac{\alpha(a+b)}{2\sqrt{ab}}+\frac{\alpha\theta}{2\sqrt{ab}}-\frac{\theta(a+b)}{4}\\
 &=-\frac38(b-a)^2-\frac{a\theta}{2}-\alpha\sqrt{\frac{a}{b}}=-\frac12\psi_\alpha(b, a)-\alpha-1=-\alpha-1
 \end{aligned}$$
 It follows
 $$U^{\nu^a_\alpha}(x)+\log x=C+A+o(\frac1x).$$
 Using the fact that $$\lim_{x\to+\infty}U^{\nu^a_\alpha}(x)+\log x=0,$$
one gets
 $$\begin{aligned}C&=-A=-\frac1{16}\Big(a^2+6ab+b^2+2(b-a)^2\log\frac{b-a}{4}\Big)\\
 &+\frac{(a+b)^2}{4}-\frac{b^2-a^2}{4}\log\frac{b-a}{4}+
 \frac{\alpha}{\sqrt{ab}}\Big(-\frac{b+a}{2}+\frac{b-a}{2}\log\frac{b-a}{4}-C_\mu\Big)
.\end{aligned}$$
 A simple computation gives the result.\\
 Moreover
  $$m_2(\nu^a_\alpha)=\int_\Bbb Rx^2\nu^a_\alpha(dx)=\frac{b-a}{128} \big(15 a^3 + 27 a^2 b +13ab^2+ 9 b^3 -16\alpha\sqrt{\frac ab}\big(3a+b\big)\big).$$
  Hence
  the energy is given by
  $$E^*_{\alpha, \sigma}=C+\frac12m_2(\nu^a_\alpha)-\alpha\int_\Bbb R\log x\,\nu^a_\alpha(dx).$$
Using mathematica we obtain
 \begin{proposition}---
 \begin{description}
 \item[(1)]For $\alpha>0$, $\sigma\geq a_c$
 $$\begin{aligned}E^*_{\alpha, \sigma}&=-\frac{1}{16} \left(a^2+6 a b+b^2+2 (a-b)^2 \log\frac{b-a}{4}-4(b+a)^2+4(b^2-a^2)\log\frac{b-a}{4}\right)
 \\&+\frac{b-a}{256} \Big(15 a^3+27 a^2 b+13 a b^2+9 b^3\Big)-\alpha\varphi(x)-\alpha^2\psi(x).\end{aligned}$$
 Where $x=\frac{b}{b-a}$, $a=\sigma$ and $b>a$ is as in proposition 2.1,$$\begin{aligned}\varphi(x)&=\frac{(b-a)^2}{8}(4x-3)\sqrt{\frac{x-1}{x}}
  +\frac12\log\Big(\frac{\sqrt x+\sqrt{x -1}}{\sqrt x-\sqrt{x-1}}\Big)\\
&+\frac{1}{2\sqrt{x(x-1)}}\big(\frac32-3x+(3-2x)\log\frac{b-a}{4}\big)  \\&+\frac{1}{4} (b-a) \Big(2-2 x+2\sqrt{(x-1)x}+2\log\Big(1+\sqrt{\frac{x-1}{x}}\Big)+\log x-(1+\log4)\Big).\end{aligned}$$
 $$\begin{aligned}\psi(x)&=\frac{-2 \big(x+\sqrt{x(x-1) }\big) \log\Big(\frac{\sqrt x+\sqrt{x-1}}{\sqrt x-\sqrt{x-1}}\Big)+\sqrt{x(x-1) } \log 4+x \log(4x(x-1))}{2 x}\\&+ \Big(\frac12-\frac12\sqrt{\frac{ x-1}{x}} \Big) \log(b-a).\end{aligned}$$
 \item[(2)]For $\alpha>0$, and all $\sigma\in[0,a_c]$,
 $$\begin{aligned}E^*_{\alpha, a_c}:=E^*_{\alpha, \sigma}&=-\frac{1}{16} \Big(3a_c^2+8 a_c b_c+3b_c^2-4a_cb_c \log\frac{b_c-a_c}{4}\big)\Big)
 \\&+\frac{b_c-a_c}{256} \Big(15 a_c^3+27 a_c^2 b_c+13 a_c b_c^2+9 b_c^3\Big)-\alpha\varphi(x_c)-\alpha^2\psi(x_c).\end{aligned}$$
 where $\displaystyle b_c=\frac23\Big(\sqrt{6(1+\alpha)-2a_c^2}-\frac{a_c}2\Big)$, $\displaystyle x_c=\frac{b_c}{b_c-a_c}$,
 and $$\begin{aligned}&\varphi(x_c)=\\&\frac1{16}(1-4\log2)+\frac{(b_c-a_c)^2}{8}(4x_c-3)\sqrt{\frac{x_c-1}{x_c}}+\frac12\log\frac{\sqrt x_c+\sqrt{x_c-1}}{\sqrt x_c-\sqrt{x_c-1}}
  \\&+\frac{1}{4} (b_c-a_c)^2\Big(\frac12\log(b_c-a_c)-2x_c+2 x^2_c+(1-2x_c)\sqrt{(x_c-1)x_c}+2\log\big(\sqrt{x_c}+\sqrt{x_c-1}\big)\Big).\end{aligned}$$
 \item[(3)]For $\alpha=0$,\\
 \begin{description}
 \item[(a)]if $\sigma\geq-\sqrt 2$, hence $a=\sigma$, $\displaystyle b=\frac23\Big(\sqrt{\sigma^2+6}+\frac \sigma2\Big),$ and
 $$E^*_{0, \sigma}=
 \frac{1}{108} \Big(81+72 \sigma^2-2 \sigma^4+(30 \sigma +2 \sigma^3) \sqrt{6+\sigma^2}-108 \log\Big(\frac{1}{6} \big(-\sigma+\sqrt{6+\sigma^2}\big)\Big)\Big).$$
 Moreover for $\sigma=0$,
 $$E^*_{0, 0}=\frac34+\frac12\log2+\frac12\log3.$$
 \item[(b)] If $\sigma\leq-\sqrt 2$,  $$E^*_{0, \sigma}=E^*_{0, -\infty}=\frac34+\frac12\log2.$$
\end{description}
\end{description}
\end{proposition}

$E^*_{0, -\infty}$ it means the energy of the equilibrium measure associate to the potential $Q_0(x)=x^2$ on $]-\infty, +\infty[$. In such a case the equilibrium measure is the semicircle law. Observe that $E^*_{0, \sigma}=E^*_{0, -\sqrt 2}$, for all $\sigma\leq-\sqrt 2$.\\
$E^*_{0, -\infty}$ is the energy of the Gaussian unitary ensemble of the eigenvalues on all the real line.\\
$E^*_{0, 0}$ can be explained as the energy of the Gaussian unitary ensemble for which all eigenvalues describe the positive real axis.\\
\section{Density of positive eigenvalues of the generalized Gaussian unitary ensemble}
We consider the generalized Gaussian unitary ensembles of random
matrices with Dyson index $\beta = 1,2,4$, corresponding
to real, complex, and quaternion entries, respectively. The
probability distribution of the entries is given by
$$\Bbb{P}_{n, \mu}(dX)=\frac1{C_n}|\det(X)|^{\beta\mu}\exp\Big({-\frac{\beta}{2}\big<X,X\big>}\Big)dX.$$
where $C_n$ is a normalizing constant, $\mu$ a postive real number and $dX$ is the Lebesgue measure on the space $H_n=Herm(n,\Bbb F)$ of hermitian matrices with respectively real, complex or quaternion coefficients $\Bbb F=\Bbb R$, $\Bbb C$, or $\Bbb H$. Consequently the joint probability density of eigenvalues
is defined on $\Bbb R^n$ by
$$\Bbb{P}_{n,\mu}(d\lambda_1,...,d\lambda_n)=\frac1{C_n}\prod_{i=1}^n|\lambda_i|^{\beta\mu}e^{-\frac{\beta}{2}
\sum\limits_{i=1}^n\lambda_i^2}|\Delta(\lambda)|^{\beta}d\lambda_1...d\lambda_n,$$
where $\Delta(\lambda)=\prod_{i<j}(\lambda_i-\lambda_j)$ is the Vandermonde determinant.

More generally one can consider this probability with an arbitrary real number $\beta>0$, but then it is not related to a matrix ensemble.

 For  $\mu\geq 0$, $\sigma>0$, consider $$\Omega_{n, \sigma}=\{x\in H_n\mid\,\lambda_{\rm min}(x)\geq \sigma\},$$ the subset of Hermitian matrices for which all it eigenvalues are bigger then some positif real number $\sigma$.  If $\mu=0$, it will be assumed that $\sigma\in\Bbb R$.

 We wish to study $\Bbb{P}_{n,\mu}(\Omega_{n, \sigma})$, the probability for a matrix $x\in H_n$, to have all it eigenvalues in $\Omega_{n, \sigma}$. It is the probability that all the eigenvalues are $\geq \sigma$.

 It can be seen by an easy translation that
 $$\Bbb{P}_{n,\mu}(\Omega_{n, \sigma})=\int_{[\sigma,+\infty[^n}\Bbb{P}_{n,\mu}(d\lambda)=\frac1{C_n}\int_{\Bbb R^n_+}\prod_{i=1}^n(\lambda_i+\sigma)^{\beta\mu}e^{-\frac{\beta}{2}
\sum\limits_{i=1}^n(\lambda_i+\sigma)^2}|\Delta(\lambda)|^{\beta}d\lambda_1...d\lambda_n.$$
Hence $$\Bbb{P}_{n,\mu}(\Omega_{n, \sigma})=\Bbb{P}^\sigma_{n, \mu}(\Omega_n),$$
where $\Omega_n$ is the cone of positive definite hermitian matrices, and $\Bbb{P}^\sigma_{n, \mu}$ is the probability on $\Bbb R^n_+$ with density
$$\Bbb{P}^\sigma_{n, \mu}(dx)=\frac1{C_n}\prod_{i=1}^n(x_i+\sigma)^{\beta\mu}e^{-\frac{\beta}{2}
\sum\limits_{i=1}^n(x_i+\sigma)^2}|\Delta(x)|^{\beta}dx_1...dx_n.$$

So the question now is to study the probability that a hermitian matrix have all its eigenvalues positive.

Let $\nu^\sigma_{n, \mu}$ be the probability measure defined on $\Bbb R_+$ by : For all bounded continuous functions $f$
$$\int_{\Bbb R_+}f(x)\nu^\sigma_{n, \mu}(dx)=\int_{\Bbb R^n_+}\frac1n\sum_{i=1}^nf(\lambda_i)\Bbb{P}^\sigma_{n, \mu}(d\lambda_1,d\lambda_2,...,d\lambda_n),$$
which means that $$\nu^\sigma_{n, \mu}=\Bbb{E}^\sigma_{n, \mu}\big(\frac1n\sum_{i=1}^n\delta_{\lambda_i}\big),$$
where $\Bbb{E}^\sigma_{n, \mu}$ is the expectation with respect the measure $\Bbb{P}^\sigma_{ n, \mu}$.

As $n$ goes to infinity we will prove after scaling the measure $\nu^\sigma_{n, \mu_n}$ by $\frac1{\sqrt n}$, it converge to some probability measure $\nu^\sigma_\alpha$, which is the statistical density of the eigenvalues bigger then $\sigma$ where $(\mu_n)_n$ is some appropriate sequence.
The means result of this section is the present theorem
\begin{theorem}---Let $\mu=(\mu_n)$ be a positive real sequence, such that $$\lim_{n\to\infty}\frac{\mu_n}{n}=\alpha.$$ Then there exist a unique $a=a(\sigma,\alpha)$, and $b=b(\sigma, \alpha)$, such that, after scaling the measure $\nu^\sigma_{n, \mu_n}$ by $\frac1{\sqrt n}$, it converge for the tight topology to the probability measure $\nu^\sigma_\alpha$ of theorem 2.4, with support $]a,b]$, and density
$f_{\alpha, \sigma}$.
This means for all bounded continuous functions $\varphi$ on $\Bbb R_+$,
$$\lim_{n\to\infty}\int_{\Bbb R_+}\varphi(\frac{x}{\sqrt n})\nu_{n,\mu_n}^\sigma(dx)=\int_{\Bbb R_+}\varphi(x)\nu^\sigma_\alpha(dx).$$
\end{theorem}
Observe that when $\alpha=0$, if we put $L(a)=b-a$, we obtain the following equation for $L(a)$,
$$\frac{3}{4}L(a)^2+aL(a)-2=0,$$
which gives $$L(a)=\frac23(\sqrt{a^2+6}-a),$$
and the density in the interval $[0, L(a)]$ is given by
$$f_{0}(x)=\frac1{2\pi}\sqrt{\frac{b-x}{x}}(2x+L(a)+2a),$$
one recover's Dean-Majumdar theorem, see for instance \cite{S}.

Furthermore when $\alpha=0$ and $a=0$, hence $b=\frac23\sqrt6$, we recover's the Dean-Majumdar density \cite{S}: the density of positive eigenvalues of a hermitian random matrix.
$$f(x)=\frac1{2\pi}\sqrt{\frac{b-x}{x}}(2x+b).$$
Note that for $\alpha>0$, if we move the barrier $\sigma$ from the right to the left, for $\sigma>a_c$, the density is given as in the theorem, moreover for $\sigma=a_c$, the density $f_{\alpha}$ vanished at the end of the support $[a_c, b_c]$ and it is the optimal density of the eigenvalues. For $0<\sigma\leq a_c$ the density keep unchanged and is given by \begin{equation}\label{022}f_{\alpha}(x)=\frac1{\pi}\sqrt{(b_c-x)(x-a_c)}\Big(1+\frac{\alpha}{\sqrt{a_cb_c}}\frac1x\Big).\end{equation} Which means that the statistical distribution of eigenvalues on all the positive real line of generalized Gaussian hermitian matrix when the dimensional of the matrix is big enough is given by $f_{\alpha}$.
Observe that as $\alpha\to 0$, since from proposition 2.1 one have, $\displaystyle \frac{b_c(\alpha)+a_c(\alpha)}{2}=\frac{\alpha}{\sqrt{b_c(\alpha)a_c(\alpha)}},$ and $0<a_c(\alpha)<\sqrt\alpha$, hence $a_c(\alpha)\to 0$, and
$$\lim_{\alpha\to 0}\frac{\alpha}{\sqrt{b_c(\alpha)a_c(\alpha)}}=\frac{b(0)}{2}=\frac23\sqrt6.$$ From equation (\ref{022}), one recover's the Dean-majumdar density of positive eigenvalues of Gaussian unitary ensemble.

\subsection{Large deviation to the left and right of $\lambda_{\min}$}
\subsection{ Case of $\alpha >0$.}
Let $$\Bbb{P}_{n,\mu_n}(\lambda_{\min}\geq t)=\frac{Z_n(t)}{Z_n(0)},$$
where $$Z_n(t)=\int_t^{+\infty}\cdots \int_t^{+\infty}e^{-\frac\beta2F_n(\lambda)}d\lambda_1...d\lambda_n,$$
and $$F_n(\lambda)=
\sum\limits_{i=1}^n\Big(\lambda_i^2+2\mu_n\log\frac{1}{\lambda_i}\Big)+2\sum_{i<j}\log\frac1{|\lambda_i-\lambda_j|}.$$
From the previous theorem, one can find the asymptotic of the density of minimal eigenvalues. In fact
by the saddle point method one obtain
$$\Bbb{P}_{n,\mu_n}\big(\lambda_{\min}(n)\geq\sqrt n\sigma\big)\approx \exp\Big(-\frac\beta2 n^2\big(E^*_{\alpha, \sigma}-E^*_{\alpha, a_c}\big)\Big)\approx \exp\Big(-\frac\beta 2n^2\Phi_+(\sigma)\Big),$$
the function $\Phi_+$ is called the right rate function.

The saddle point method is not able to capture
fluctuations to the left of $\lambda_{\min}$. Since we only consider leading terms $O(n^2)$, which
capture bulk properties. the asymptotic density $f_{\alpha, \sigma}$ sea is a priori not subject
to forces capable of macroscopic rearrangements. Following
this physical picture, the left rate function is determined by
the energy cost in pulling the leftmost charge in the external
potential of the Coulomb gas $F_n(\lambda)$ and the interaction of the
charge with the unperturbed $f_{\alpha, \sigma}$ sea. This energy cost for
$\lambda_{\min}=t \ll a_c\sqrt n$ can be estimated for large $n$ using $F_n(\lambda)$
$$\Delta E_n(t)=F_n(t,\lambda_2,...,\lambda_n)-F_n(a_c\sqrt n,\lambda_2,...,\lambda_n).$$
Hence$$\begin{aligned}\Delta E_n(t)&=t^2-2\mu_n\log t-2\sum_{j=1}^n\log{|t-\lambda_j|}+ C_n\\&
=t^2-2\mu_n\log t-2n\int\log|t-\lambda|\nu_{n,\mu_n}^\sigma(d\lambda)+C_n.\end{aligned}$$
where $C_n$ is given by the condition $\displaystyle\Delta E(t=a_c\sqrt n)=0$ and $0<\sigma<a_c$.

For $n$ large enough using the fact that $\mu_n\sim\alpha n$, and the convergence of the measure $\nu_{n,\mu_n}^\sigma$ tightly to $f_\alpha$, one gets for the new scaling variable, $t=x\sqrt n$,
$$\frac{\Delta E_n(t)}{n} \sim \Delta E_\alpha(x)=x^2-2\alpha\log x-2\int_{a_c}^{b_c}\log|x-\lambda|f_\alpha(\lambda)d\lambda+C,$$
where $$f_{\alpha}(\lambda)=\frac1{2\pi}\sqrt{(b_c-\lambda)(\lambda-a_c)}(1+\frac{\alpha}{\sqrt{a_cb_c}}\frac1\lambda).$$
Hence $$\Delta E_\alpha(x)=x^2-2\alpha\log(x)+2U^{\nu_\alpha^{a_c}}(x)+C,$$
Such expression coincide with the previous one of  proposition 3.2. Hence
$$\Delta E_\alpha(x)=2\int_x^{a_c}h_{\alpha, a_c}(\lambda)d\lambda,$$
where $$h_{\alpha, a_c}(\lambda)=\sqrt{(b_c-\lambda)(a_c-\lambda)}\Big(1+\frac{\alpha}{\sqrt{a_cb_c}}\frac1\lambda\Big).$$
A simple computation gives
$$\begin{aligned}\Delta E_\alpha(x)&=\frac{1}{4 }(a_c-b_c)^2 \Big(-\xi \sqrt{-1+\xi^2}-\log\big(-\xi+\sqrt{-1+\xi^2}\big)\Big)\\&+\frac{1}{2}(a^2_c-b^2_c)  \Big(\sqrt{-1+\xi^2}-\frac{3}{2} \sqrt{-1+r^2} \log\big(-1+r^2\big)+r \log\big(-\xi+ \sqrt{-1+\xi^2}\big)\Big)\\&-\frac{1}{2}(a^2_c-b^2_c) \sqrt{r^2-1} \log\Big(\frac{-1-r\xi+\sqrt{\big(-1+r^2\big) \big(-1+\xi^2\big)}}{(r+\xi)\big(-1+r^2\big)^{3/2}} \Big),\end{aligned}$$
where $\displaystyle\xi=\frac{2x-b_c-a_c}{b_c-a_c}$, $ \displaystyle r=\frac{b_c+a_c}{b_c-a_c}$, and $x\in]0,a_c[$.
Moreover for the probability,
one gets
$${\cal P}_{n, \mu_n}(\lambda_{\min}(n)\leq x\sqrt n)\sim e^{-\frac{n\beta}2\Delta E_\alpha(x)},$$
Since $
x\sim a_c$, which means that $\xi\sim -1$, hence

$$\Delta E_\alpha(x)=-\frac{(a_c-b_c)^2}{4} \left(\frac{4  \sqrt{2}}{3}+\frac{4 \left(-\frac{1}{2 \sqrt{2}}+\frac{ r}{6 \sqrt{2}}+\frac{(1+r) \left(-5 \sqrt{2}+\sqrt{2} r\right) }{12 (-1+r) }\right) \alpha }{(b_c-a_c) \sqrt{a_c b_c}}\right) (-\xi-1)^{3/2}+O[-\xi-1]^{5/2}.$$
for $\xi\leq -1$. Substituting the expression of $\xi$, and $x$ in the previous formula, it yields
for $\frac{t-a_c\sqrt n}{\sqrt n}\ll1 $,
$${\cal P}_{n, \mu_n}\big(\lambda_{\min}(n)\leq t\big)\sim e^{-\frac{n\beta}2C_0\Big(\frac{a_c\sqrt n-t}{\sqrt n}\Big)^{\frac 32}},$$

\begin{equation}\label{pq}{\cal P}_{n, \mu_n}\big(\lambda_{\min}(n)\leq t\big)\sim \exp\Big({-\frac{\beta}2C_0\Big(n^{\frac16}(a_c\sqrt n-t)\Big)^{\frac 32}}\Big),\end{equation}
where $$C_0:=C_0(\alpha)=\sqrt{b_c-a_c} \bigg(\frac{4}{3}+\frac{\Big(-1+\frac{ r}{3}+\frac{(1+r) (-5 +r) }{3 (-1+r) }\Big) \alpha }{(b_c-a_c) \sqrt{a_c b_c}}\bigg).$$
Which can be simplified as
$$C_0=\sqrt{b_c-a_c} \bigg(\frac{4}{3}+\frac{r\big(r-2-\sqrt 5\big)\big(r-2+\sqrt 5\big)}{6(r-1)}\bigg),\quad r=r(\alpha)=\frac{b_c+a_c}{b_c-a_c} .$$
Such result coincide exactly with the left tail behavior of the Tracy-Widom limiting distribution.
In fact if we put $\lambda_{\min}(n)=a_c\sqrt n-\frac34C_0^{\frac23}n^{-\frac16}\chi$, the $\chi$ has an $n$-independent PDF $g_\beta$, from the left of $a_c$
and $$g_\beta(\chi)=\exp\big(-\frac{2\beta}3\chi^{\frac 32}\big)\quad {\rm for }\;\;\chi\rightarrow+\infty.$$
\subsection{ Case of $\alpha=0$}
 If $\mu_n=0$, we are in the case of the Gaussian unitary ensemble which is treated by Dean-Majumdar \cite{S}, where they studies the rate function for the maximal eigenvalue.  Here we assume more general setting where we suppose that $\displaystyle\lim_{n\to+\infty}\frac{\mu_n}{n}=\alpha=0$. we'll see that the result of Dean-Majumdar still unchanged. We will gives the rate function for the minimal eigenvalue.
In such a case we saw that the limit density of eigenvalues which are in the interval $[\sigma,+\infty[$, for all $\sigma\in\Bbb R$, is
$$f_0(x)=\frac1{2\pi}\sqrt{\frac{b-x}{x-\sigma}}\Big(2x+b-\sigma\Big),$$
for all $\sigma\in\Bbb R$, where $b=b(\sigma)=\frac23\big(\sqrt{\sigma^2+6}+\frac\sigma2\big)$.
One gets the probability asymptotic
$$\Bbb{P}_{n,\mu_n}\big(\lambda_{\min}(n)\geq\sqrt n\sigma\big)\approx \exp\Big(-\frac\beta2 n^2\big(E^*_{0, \sigma}-E^*_{0}\big)\Big)\approx \exp\Big(-\frac\beta 2n^2\Phi_+(\sigma)\Big),$$
Here $E^*_0$ is the full energy in the real line. It is given by
$$E^*_0=E^*_{0,-\infty}=E^*_{0,-\sqrt 2}=\frac34+\frac12\log 2.$$
 In the same way from  proposition 3.5  one has the expression of $E^*_{0,\sigma}$,
$$ E^*_{0,\sigma}=\frac{1}{108} \Big(81+72 \sigma^2-2 \sigma^4+\big(30\sigma +2 \sigma^3\big) \sqrt{6+\sigma^2}-108 \log \Big(\frac{1}{6} \big(-\sigma+\sqrt{6+\sigma^2}\big)\Big)\Big).$$

 Hence we obtain
 $$\Bbb{P}_{n,\mu_n}\big(\lambda_{\min}(n)\geq\sqrt n\sigma\big)\approx \exp\Big(-\frac\beta 2n^2\Phi_+(\sigma)\Big),$$
 where
 $$\Phi_+(\sigma)=\frac1{54}\Big(36\sigma^2-\sigma^4+\big(15\sigma+\sigma^3\big)\sqrt{\sigma^2+6}+
 27\big(\log18-2\log\big(\sqrt{\sigma^2+6}-\sigma\big)\Big).$$
For $\sigma\rightarrow (-\sqrt 2)^+$ one gets
 $$\Phi_+(\sigma)\sim \frac{\sqrt 2}{6}(\sigma+\sqrt 2)^3.$$
Thus for the scaling variable $t$, $\displaystyle t=\sigma\sqrt n$, and $\frac{t+\sqrt{2n}}{\sqrt n}\ll 1$, $t\geq-\sqrt{2 n}$ we obtain
 $$\Bbb{P}_{n,\mu_n}\big(\lambda_{\min}(n)\geq t\big)\sim \exp\Big(-\frac{\beta\sqrt 2}{12}\big(n^{\frac 16}(t+\sqrt {2n})\big)^3\Big).$$
 Note that this matches exactly with the right tail behavior of the Tracy-Widom limiting distribution. The same hold for $\alpha>0$, but the expression of $E_{\alpha, \sigma}$ is more complicate.

 To see the left rate function, we have for $\sigma\leq -\sqrt 2$, from the previous case with ($a_c=-\sqrt 2, b_c=\sqrt 2$), that
$$\Phi_-(\sigma)=\Delta E_0(\sigma)= \log 2-\sigma \sqrt{\sigma^2-2}-2 \log\big(-\sigma+\sqrt{\sigma^2-2}\big).$$
 For $\sigma\rightarrow (-\sqrt 2)^-$ one gets
 $$\Phi_-(\sigma)\sim \frac{2^{\frac {11}4}}{3}(-\sigma-\sqrt 2)^{\frac32}.$$
  For $\sigma\mapsto -\infty$
 $$\Phi_-(\sigma)\sim \sigma^2.$$
 Since
 $$\Bbb{P}_{n,\mu_n}\big(\lambda_{\min}(n)\leq \sigma\sqrt n\big)\sim\exp\big(-n\frac{\beta}{2}\Phi_-(\sigma)\big),$$
 and
  This means with the new scale variable $t=\sigma\sqrt n\leq-\sqrt{2n}$ that,
   $$\Bbb{P}_{n,\mu_n}\big(\lambda_{\min}(n)\leq t\big)\sim \exp\Big(-\frac{\beta 2^{\frac 74}}{3}\big(n^{\frac 16}\big|t+\sqrt {2n}\big|\big)^{\frac32}\Big).$$
  This matches exactly with the left tail behavior of the Tracy-Widom limiting distribution.
 Such result can be deduced from equation (\ref{pq}). In fact for $\alpha=0$, hence $a_c(\alpha)=-\sqrt 2$, $b_c(\alpha)=\sqrt 2$ and $C_0(\alpha)=\frac{2^{\frac {11}4}}{3}$.
 \subsection{Proof of theorem 5.1}
 In the rest of the section we will prove theorem 4.1, for this purpose we need some preliminary results.

For a positive real sequence $\mu=(\mu_n)_n$ and $\sigma$ as in the previous, let $K^\sigma_n$ be the function on $\Bbb R_+^n$, defined by
$$K^\sigma_n(x)=\sum_{i\neq j}\log\frac{1}{|x_i-x_j|}+(n-1)\sum_{i=1}^nQ^\sigma_{\alpha_n}(x_i),$$
where $\displaystyle Q^\sigma_{\alpha_n}(x)=(x+\sigma)^2+2\alpha_n\log\frac1{x+\sigma},$ $\displaystyle\alpha_n=\frac{\mu_n}{n},$ and
$$k^\sigma_n(x,y)=\log\frac{1}{|x-y|}+\frac12Q^\sigma_{\alpha_n}(x)+\frac12Q^\sigma_{\alpha_n}(y).$$
The function $K^\sigma_n$ is bounded from below, moreover $\lim_{x\to+\infty}K^\sigma_n(x)=+\infty$, it follows that $K^\sigma_n$ attaint it minimum at a point say, $x^{(n, \sigma)}=(x_1^{(n, \sigma)},...,x_n^{(n, \sigma)} )$.
Let $$\tau^\sigma_n=\frac1{n(n-1)}\inf_{\Bbb R_+^n} K^\sigma_n(x),$$
and $$\rho^\sigma_n=\frac1n\sum_{i=1}^n\delta_{x^{(n, \sigma)}_i}.$$

For a probability measure $\mu$ on $\Bbb R_+$, and $\delta\geq 0$, one consider the energy
$$E_{\delta, \sigma}(\mu)=\int_{\Bbb R_+^2}\log\frac{1}{|s-t|}\mu(ds)\mu(dt)+\int_{\Bbb R_+}Q_\delta^\sigma(s)\mu(ds).$$
and $$E^*_{\delta, \sigma}=\inf_{\mu}E_{\delta, \sigma}(\mu),$$
where the minimum is taken over all compactly support measures with support in $\Bbb R_+$.
Moreover, defined the scaled density
 $$\Bbb{P}^\sigma_{ n, \mu_n}(dx)=\frac1{Z_n}\prod_{i=1}^n(x_i+\sigma)^{\beta\mu_n}e^{-n\frac{\beta}{2}
\sum\limits_{i=1}^n(x_i+\sigma)^2}|\Delta(x)|^{\beta}dx_1...dx_n.$$
where $Z_n$ is a normalizing constant.
\begin{proposition}---
Let $\mu=(\mu_n)$ be a positive real sequence, and $\sigma>0$,  such that $$\lim_{n\to\infty}\frac{\mu_n}{n}=\alpha.$$  If $\alpha=0$, it will be assumed that $\sigma\in\Bbb R$. Then there exist $a=a(\sigma,\alpha)$, $(a\geq\sigma)$, and $b=b(\sigma, \alpha)$, $(b>a)$, such that \begin{description}
\item[(1)] $\displaystyle\lim_{n\to\infty}\tau^\sigma_n=E^*_{\alpha,\sigma}=E_{\alpha, \sigma}(\nu^a_\alpha).$
\item[(2)] The measure $\rho^\sigma_n$ converge for the tight topology to the equilibrium measure $\nu^a_\alpha$.
\item[(3)] $\displaystyle\lim_{n\to\infty}-\frac1{n^2}\log Z_n=\frac\beta2E^*_{\alpha, \sigma}$.
\end{description}
\end{proposition}
\begin{lemma} For every $n\in\Bbb N$, let $\displaystyle\alpha_n=\frac{\mu_n}{n}$ and assume $\lim\limits_{n\to+\infty}\alpha_n=\alpha$. Then there exist $a_n\geq \sigma$, $b_n>a_n$, such that
\begin{description}
\item[(1)] $E^*_{\alpha_n, \sigma}=E_{\alpha_n, \sigma}(\nu_{\alpha_n}^{a_n})$, $\lim\limits_{n\to\infty }a_n=a$, and   $\lim\limits_{n\to\infty}b_n=b$.
\item[(2)] $\lim\limits_{n\to\infty}E^*_{\alpha_n, \sigma}=E^*_{\alpha, \sigma}$.
\end{description}
\end{lemma}
{\bf Proof.}---\\
{\bf Step (1).} From proposition 1.1 and 3.1, it follows that with $\Sigma=[\sigma,+\infty[$, and $Q_{\alpha_n}$, there is $a_n\geq\sigma$ and $b_n>a_n$, such that
the measure which realize the minimum of the energy $E_{\alpha_n, \sigma}$ is $\nu_{\alpha_n}^{a_n}$ with support $]a_n, b_n]$, where $a_n$ and $b_n$ are the unique solution of the two equations $\psi_{\alpha_n}(b_n, a_n)=0$ and $\varphi_{\alpha_n}(b_n,a_n)\geq 0$ on $[\sigma,+\infty[$.
Hence
$$\frac34(b_n-a)^2+a(b_n-a)+2\alpha_n\frac{\sqrt a}{\sqrt b_n}-2\alpha_n-2=0,\;a_n+b_n-\frac{2\alpha_n}{\sqrt{a_n b_n}}\geq 0$$
it follows $$\sigma\leq a_n\leq b_n\leq\frac43(2\alpha_n+2),$$
Since the sequence $\alpha_n$ converge to $\alpha$, hence $a_n$ and $b_n$ are bounded and there is some subsequence $a_{n_k}$, $b_{n_k}$ convergent respectively to $a_0$ and $b_0$ and $b_0\geq a_0\geq \sigma$. It follows that
$$\frac34(b_0-a)^2+a(b_0-a)+2\alpha\frac{\sqrt a}{\sqrt b_0}-2\alpha-2=0,\;a_0+b_0-\frac{2\alpha}{\sqrt{a_0 b_0}}\geq 0$$
which means that $$\psi_\alpha(b_0, a_0)=0,\;{\rm and},\;\varphi_\alpha(b_0, a_0)\geq 0.$$
Moreover we saw that the solution $(x,y)$ of the equations $\psi_\alpha(y, x)=0$ and $\varphi_\alpha(y, x)\geq 0$ in the interval $[\sigma, +\infty[$ is unique the unique pair  $(a,b)$. It follows that
$a_0=a$ and $b_0=b$, and $a$, $b$ are the unique limits for a subsequences of $a_n$ respectively $b_n$, hence these sequences $a_n$ and $b_n$  converge to $a$ respectively $b$.\\
{\bf step (2).} We saw that supp$(\nu_{\alpha_n}^{a_n})\subset[\sigma,+\infty[$. Hence
$$E^*_{\alpha, \sigma}=E_{\alpha, \sigma}(\nu^a_{\alpha})\leq E_{\alpha, \sigma}(\nu^{a_n}_{\alpha_n}),$$
furthermore, we have seen in {\bf step(1)} that $E^*_{\alpha_n, \sigma}=E_{\alpha_n, \sigma}(\nu_{\alpha_n}^{a_n})$, hence
$$E_{\alpha, \sigma}(\nu^{a_n}_{\alpha_n})=E^*_{\alpha_n, \sigma}+\int_{\Bbb R_+}\Big(Q^\sigma_\alpha(x)-Q^{\sigma}_{\alpha_n}(x)\Big)\nu^{a_n}_{\alpha_n}(dx),$$
By using the dominate convergence theorem and the fact that $a_n$, $b_n$ converge, we deduce that
$$\lim_{n\to\infty}\int_{\Bbb R_+}\Big(Q^\sigma_\alpha(x)-Q^{\sigma}_{\alpha_n}(x)\Big)\nu^{a_n}_{\alpha_n}(dx)=0.$$
and \begin{equation}\label{p}E^*_{\alpha, \sigma}\leq\liminf_nE^*_{\alpha_n, \sigma}.\end{equation}
Furthermore we know that supp$(\nu_\alpha^a)\subset[\sigma,+\infty[$, it follows that
\begin{equation}\label{q}E^*_{\alpha_n, \sigma}\leq E_{\alpha_n, \sigma}(\nu_\alpha^a).\end{equation}
Moreover
 $$E_{\alpha_n, \sigma}(\nu^a_\alpha)=E_{\alpha, \sigma}(\nu_\alpha^a)+\int_\Bbb R\Big(Q^\sigma_{\alpha_n}(x)-Q^\sigma_{\alpha}(x)\Big)\nu^a_{\alpha}(dx),
.$$
 Hence
 $$E_{\alpha_n, \sigma}(\nu^a_\alpha)=E_{\alpha, \sigma}^*+\int_\Bbb R\Big(Q^\sigma_{\alpha_n}(x)-Q^\sigma_{\alpha}(x)\Big)\nu^a_{\alpha}(dx),$$
 by the dominated convergence theorem, one gets
 $$\lim_{n\to +\infty}E_{\alpha_n, \sigma}(\nu^a_\alpha)=E_{\alpha, \sigma}^*,$$
 and from (\ref{q}) we obtain \begin{equation}\label{r}\limsup_nE^*_{\alpha_n, \sigma}\leq E_{\alpha, \sigma}^*.\end{equation}
 Thus from equations (\ref{p}) and (\ref{r}), it follows that
 $$E^*_{\alpha, \sigma}\leq\liminf_nE^*_{\alpha_n, \sigma}\leq\limsup_nE^*_{\alpha_n, \sigma}\leq E_{\alpha, \sigma}^*.$$
 Which gives the desired result.\\
{\bf Proof of proposition 5.2.}---\\
{\bf Step 1 and 2:}
For a probability measure $\mu$,
$$\int_{\Bbb R^n_+}K^\sigma_n(x)\mu(dx_1)...\mu(dx_n)=n(n-1)\int_{\Bbb R^2_+}\log\frac{1}{|x-y|}\mu(dx)\mu(dx)+n(n-1)\int_{\Bbb R_+}Q^\sigma_{\alpha_n}(x)\mu(dx),$$
hence $$\tau_n^\sigma\leq E_{\alpha_n, \sigma}(\mu).$$
For $\mu=\nu^{a_n}_{\alpha_n}$,
 \begin{equation}\label{34}\tau_n^\sigma\leq E^*_{\alpha_n, \sigma}.\end{equation}
Moreover $$K^\sigma_n(x^{(n, \sigma)})=\sum_{i\neq j}k^\sigma_n\Big(x_i^{(n, \sigma)},x_j^{(n, \sigma)}\Big)\geq(n-1)\frac12\Big(\sum_{i=1}^nh^\sigma_{\alpha_n}\Big(x_i^{(n, \sigma)}\Big)+
\sum_{i=1}^nh_{\alpha_n}^\sigma\Big(x_i^{(n, \sigma)}\Big)\Big),$$
where $\displaystyle h_{\alpha_n}^\sigma(x)=Q^\sigma_{\alpha_n}(x)-\log(1+x^2),$
Since $$\int_{\Bbb R_+}h_{\alpha_n}^\sigma(t)\rho^\sigma_n(dt)=\frac1n\sum_{i=1}^nh_{\alpha_n, \sigma}\Big(x_i^{(n, \sigma)}\Big),$$
it follows that,
$$\int_{\Bbb R_+}h_{\alpha_n}^\sigma(t)\rho^\sigma_n(dt)\leq\tau^\sigma_n\leq  E^*_{\alpha_n, \sigma}.$$
Moreover, by the convergence of $\alpha_n$ to $\alpha$, there is some constant $c$ such that for all $t$
$$h_{\alpha_n}^\sigma(t)\geq h_c^\sigma(t).$$
Furthermore from {\bf Step(1)} of the previous lemma, there is some constant $C$, such that
$$\int_{\Bbb R_+}h_{c}^\sigma(t)\rho^\sigma_n(dt)\leq E^*_{\alpha_n, \sigma}\leq C,$$
using the fact that $\lim_{x\to+\infty}h_c^\sigma(x)=+\infty,$ then by the Prokhorov criterium there is some subsequence $\rho_{n_k}$, which convergent to $\rho$ for the tight topology.

For $\ell\geq 0$, let $k_n^{\sigma, \ell}(x,y)=\inf(k^\sigma_n(x,y),\ell)$,
defined $$E^\ell_{\alpha_n, \sigma}(\mu)=\int_{\Bbb R^2_+}k^{\sigma, \ell}_n(x,y)\mu(dx)\mu(dy),$$
Divide $\Bbb R_+^2$ to four regions,
$$R_1=\{(x,y)\in\Bbb R^2_+\mid x\leq 1-\sigma,\,{and}\,y\leq 1-\sigma\},$$
$$R_2=\{(x,y)\in\Bbb R^2_+\mid x\geq 1-\sigma,\,{and}\,y\geq 1-\sigma\},$$
$$R_3=\{(x,y)\in\Bbb R^2_+\mid x\leq 1-\sigma,\,{and}\,y\geq 1-\sigma\},$$
$$R_4=\{(x,y)\in\Bbb R^2_+\mid x\geq 1-\sigma,\,{and}\,y\leq 1-\sigma\},$$
Since for every $\e>0$, there is $n_0$, such that for all $n\geq n_0$,
$$\alpha-\e\leq\alpha_n\leq\alpha+\e,$$ it follows that
for $(x,y)\in R_1$,
$$k^\sigma_{\alpha-\e}(x,y)\leq k_n^\sigma(x,y),$$
for $(x,y)\in R_2$,
$$k^\sigma_{\alpha+\e}(x,y)\leq k_n^\sigma(x,y),$$
for $(x,y)\in R_3$,
$$\frac12k^\sigma_{\alpha+\e}(x,y)+\frac12k^\sigma_{\alpha-\e}(x,y)\leq k_n^\sigma(x,y),$$
the last result is valid in $R_4$ by symmetry. If we set $(\theta_1, \theta_2)=(1, 0), (0, 1), (\frac12, \frac12)$, in respectively for $R_1$, $R_2$ and $R_3\cup R_4$. It follows that
for all $(x,y)\in\Bbb R^2_+$,
$$\theta_1 k^\sigma_{\alpha-\e}(x,y)+\theta_2k^\sigma_{\alpha+\e}(x,y)\leq k_n^\sigma(x,y).$$
Take the infimum it yields
$$\theta_1 k^{\sigma,\ell}_{\alpha-\e}(x,y)+\theta_2k^{\sigma, \ell}_{\alpha+\e}(x,y)\leq k_n^{\sigma, \ell}(x,y),$$
and $$\theta_1 E^\ell_{\alpha-\e, \sigma}(\rho^\sigma_{n_k})+\theta_2E^\ell_{\alpha+\e, \sigma}(\rho^\sigma_{n_k})\leq E^\ell_{\alpha_n, \sigma}(\rho^\sigma_{n_k})\leq\tau_n^\sigma+\frac{\ell}{n}.$$
Thus from equation (\ref{34}), we have
\begin{equation}\label{35}\theta_1 E^\ell_{\alpha-\e, \sigma}(\rho^\sigma_{n_k})+\theta_2E^\ell_{\alpha+\e, \sigma}(\rho^\sigma_{n_k})\leq\tau_n^\sigma+\frac{\ell}{n}\leq E^*_{\alpha_n, \sigma}+\frac{\ell}{n}.\end{equation}
The cut kernel $k^{\sigma, \ell}_{\alpha\pm\e}(x,y)$ is bounded and continuous, and the probability measure
 $\rho^\sigma_{n_k}$ converge tightly to $\rho^\sigma$, hence
  $\displaystyle\lim_{k\to+\infty}E_{\alpha\pm\e}(\rho^\sigma_{n_k})=E_{\alpha\pm\e}(\rho^\sigma)$. Thus by the previous lemma one gets
 $$\theta_1 E^\ell_{\alpha-\e, \sigma}(\rho^\sigma)+\theta_2E^\ell_{\alpha+\e, \sigma}(\rho^\sigma)\leq E^*_{\alpha, \sigma}.$$
As $\ell$ goes to $+\infty$, by the monotone convergence theorem one obtains
$$\theta_1 E_{\alpha-\e, \sigma}(\rho^\sigma)+\theta_2E_{\alpha+\e, \sigma}(\rho^\sigma)\leq E^*_{\alpha, \sigma}.$$
thus,
$$E_{\alpha, \sigma}(\rho^\sigma)\leq E^*_{\alpha, \sigma}.$$

By the definition of the equilibrium measure we obtains $E^*_{\alpha, \sigma}=E_{\alpha, \sigma}(\rho^\sigma)=E_{\alpha, \sigma}(\nu^a_\alpha)$, it follows by unicity of the equilibrium measure that $\rho^\sigma=\nu^a_\alpha$. Thus  the only possible limit for a subsequence of $\rho_n^\sigma$ is $\nu^a_\alpha$, hence the sequence $\rho_n$ it self converge to $\nu^a_\alpha$.
Moreover from equation (\ref{35}) one gets
 $$\lim\limits_n\tau_n=E^*_{\alpha, \sigma}.$$
{\bf Step 3:}
We saw for every $x\geq 0$, $K^\sigma_n(x)\geq n(n-1)\tau^\sigma_n,$ hence
$${Z}_n\leq e^{-\frac\beta 2n(n-1)\tau_n}\Big(\int_0^{+\infty} e^{-\frac\beta 2Q^\sigma_{\alpha_n}(x)dx}\Big)^n,$$
moreover
$$\int_0^{+\infty} e^{-\frac\beta 2Q^\sigma_{\alpha_n}(x)dx}=\int_0 ^{+\infty}(x+\sigma)^{\beta\alpha_n}e^{-\frac\beta 2(x+\sigma)^2}dx\leq\frac12\big(\frac{2}{\beta}\big)^{\frac{\beta\alpha_n}{2}}\Gamma\Big(\frac{\beta\alpha_n}{2}+\frac12\Big),$$
hence $$\frac1{n^2}\log{Z}_n\leq -\frac{n-1}{n}\tau_n+\frac1n\log\Big(\frac12\big(\frac{2}{\beta}\big)^{\frac{\beta\alpha_n}{2}}\Gamma\Big(\frac{\beta\alpha_n}{2}
+\frac12\Big)\Big).$$
Here we used the fact that $\alpha_n$ converge. Then
$$\limsup_n\frac1{n^2}\log{Z}_n\leq -\frac\beta 2E^*_{\alpha, \sigma}.$$
Furthermore $${Z}_n\geq \int_{\Bbb R^n}e^{-\frac\beta 2K_n(x)-\frac\beta 2Q^\sigma_{\alpha_n}(x)-\sum\limits_{i=1}^n\log f_{\alpha}(x_i)}\prod_{i=1}^n\nu^a_\alpha(dx_i),$$
Applying Jensen's inequality we obtain

$${Z}_n\geq\exp\int_{\Bbb R^n}\Big(-\frac\beta 2K_n(x)-\frac\beta 2Q^\sigma_{\alpha_n}(x)-\sum\limits_{i=1}^n\log f_{\alpha}(x_i)\Big)\prod_{i=1}^n\nu^a_\alpha(dx_i),$$
hence
$${Z}_n\geq e^{-\frac\beta 2\big(n(n-1)E^*_{\alpha, \sigma}\big)}\exp\Big(-\frac\beta 2n\int_a^bQ^\sigma_{\alpha_n}(x)f_{\alpha}(x)dx\Big)\exp\Big(-n\int_a^bf_{\alpha}(x)\log f_{\alpha}(x)\Big)dx,$$
since $|Q^\sigma_{\alpha_n}(x)|\leq b^2+c\,\max(|\log b|, |\log a|)$, where $c$ is some bound for the sequence $\alpha_n$. Moreover the function $x\mapsto f_{\alpha}(x)\log f_{\alpha}(x)$ is continuous on $[a,b]$. Hence
$$\liminf_n\frac1{n^2}\log{Z}_n\geq -\frac\beta 2E^*_{\alpha, \sigma},$$
and the conclusion hold
 $$-\frac\beta 2E^*_{\alpha, \sigma}\leq\liminf_n\frac1{n^2}\log{Z}_n\leq\limsup_n\frac1{n^2}\log{Z}_n\leq -\frac\beta 2E^*_{\alpha, \sigma}.$$
 {\bf Proof of theorem 5.1.}---The proof of the theorem follow as the proof in (\cite{F}, Faraut, theorem IV.5.1).
 \begin{remark} Let $\alpha\geq 0$, and $\sigma<0$. If $\alpha=0$, it is assumed that $\sigma\in\Bbb R$. we restrict the probability density $\Bbb {P}_{n,\mu_n}$ to the set $\Lambda_\sigma=\Big\{x\in H_n\mid\lambda_{\max}(x)\leq\sigma\Big\}$, where $\lambda_{\max}(x)$ is the maximal eigenvalues of the hermitian matrix $x$. Under the same condition of theorem 4.1, one can prove that the measure $\nu_{n,\mu_n}^\sigma$ rescaled by $\frac{1}{\sqrt n}$ converge in the tight topology to same probability measure $\nu_\alpha^\sigma$ with density $g_\alpha$ and support $[a,b[$
 $$g_\alpha(x)=\frac1{2\pi}\sqrt{\frac{x-a}{b-x}}\Big(-2x+b-a+2\alpha\sqrt{\frac{b}{a}}\frac1x\Big).$$
 where  $b=b(\alpha, \sigma)\leq\sigma$ and $a=a(\sigma,\alpha)<b$ is the unique solution of the following equations
$$a+b+\frac{2\alpha}{\sqrt{ba}}\leq 0,$$
$$\frac{3}{4}(a-b)^2-b(b-a)+2\alpha{\sqrt\frac{  b}{ a}}-2\alpha-2=0.$$

 \end{remark}
\section{Computation of the probability of positive eigenvalues}
In this section we will compute the asymptotic of the probability $p_n(\alpha)$, for which all the eigenvalues are positive. For $\sigma=0$, we have
 $$p_n(\mu_n)=\frac{\Bbb{P}_n^0(\Omega_n)}{\Bbb P_n(H_n)},$$
Hence $$p_n(\mu_n)=\frac{Z^0_n}{Z_n},$$
where
$$Z_n^0=\int_{\Bbb R^n_+}\prod_{i=1}^nx_i^{\beta\mu_n}e^{-n\frac{\beta}{2}
\sum\limits_{i=1}^nx_i^2}|\Delta(x)|^{\beta}dx_1...dx_n.$$
and
$$Z_n=\int_{\Bbb R^n}\prod_{i=1}^nx_i^{\beta\mu_n}e^{-n\frac{\beta}{2}
\sum\limits_{i=1}^nx_i^2}|\Delta(x)|^{\beta}dx_1...dx_n.$$
It can be seen that
$$\log p_n(\mu_n)=\log Z_n^0-\log Z_n,$$
hence
$$\lim_{n\to+\infty}\frac1{n^2}\log p_n(\mu_n)=-\frac\beta2\Big(-E^*_{\alpha}+E^*_{\alpha, a_c}\Big):=\lim_{n\to+\infty}\frac1{n^2}\log p_n(\alpha).$$

{\bf Small value of $\alpha$:}
Since $$E^*_\alpha=\frac34+\frac12\log2+(\frac{3}{2}+\log2)\alpha+\alpha^2\log2\alpha-(\alpha^2+\alpha+\frac14)\log(1+2\alpha).$$
See for instance \cite{bo}.
Hence for $\alpha$ small enough $$E^*_\alpha=\frac34+\frac12\log2+(1+\log2)\alpha+o(\alpha),$$
Moreover expanding for $\alpha$ near zero, the energy $E_{\alpha,a_c}$, then we obtains, from proposition 4.5. and theorem 2.4.

 $$E_{\alpha, a_c}=\frac34  + \frac12\log6 +
 C\alpha+o(\alpha).$$
    where  $\displaystyle C=\frac1{432}\Big(-36 (-6 + \sqrt6) + (54 - 161 \sqrt6) \log2 +
    27 (10 + \sqrt6) \log3\Big)\approx0.6045$

  \begin{proposition}---\
  \begin{description}
 \item[(1)] For $\alpha>0$, $$\lim_{n\to\infty}\frac1{n^2}\log p_n(\alpha)=-\frac\beta4\log3 +\frac\beta2 (1+\log 2-C)\alpha+o(\alpha),$$
     where $o(\alpha)$ is a small terms in $\alpha$.
     \item[(2)] For $\alpha=0$,
    $$ \lim_{n\to\infty}\frac1{n^2}\log p_n(0)=-\frac\beta4\log3.$$
 \end{description}
  \end{proposition}
    Hence $p_n(\alpha)$ converge to zero very rapidly as $n\to+\infty$. Moreover for small values of $\alpha>0$, it is like $e^{-(c+o(\alpha))n^2}$ where $c=\frac\beta4\log3 +\frac\beta2 (C-1 -\log2)\alpha$. But for $\alpha=0$, $c=\frac\beta4\log3$.
    The second step of the proposition is du to Dean-Majumdar \cite{S}.

{\bf Example:} For example for $\sigma=0$, and $\alpha=0.1$, the support of the measure $\nu_1^0$ is $[a_c,b_c]=[0.00796, 1.71004]$, and $$\nu_{0.1}^0(dx)=\frac{1}{\pi}\sqrt{(1.71004-x)(x-0.00796)}\Big(1+\frac{0.858}{x}\Big).$$
Exact value of $p_n$ from energy expression:
$$E_{\alpha, a_c}=1.869, E_\alpha=1.23416, $$
$$ p_n(0.1)\approx_{n\to+\infty} e^{-(0.3174)\beta n^2}$$
Value of $p_n$ using the developmental for $\alpha$ near $0$, $(\alpha=0.1)$.
$$E_{\alpha, a_c}\approx 1.6246, E_\alpha\approx 1.2658,$$
$$ p_n(0.1)\approx_{n\to+\infty} e^{-(0.2202)\beta n^2}$$
The error is of order $0.0972$.
$$ p_n(0)\approx_{n\to+\infty} e^{-(0.2746)\beta n^2}$$
 \begin{center}{\bf Plotting of the density $f_{\alpha, \sigma}$}\end{center}

\begin{figure}[h]
\centering\scalebox{0.6}{\includegraphics[width=15cm, height=10cm]{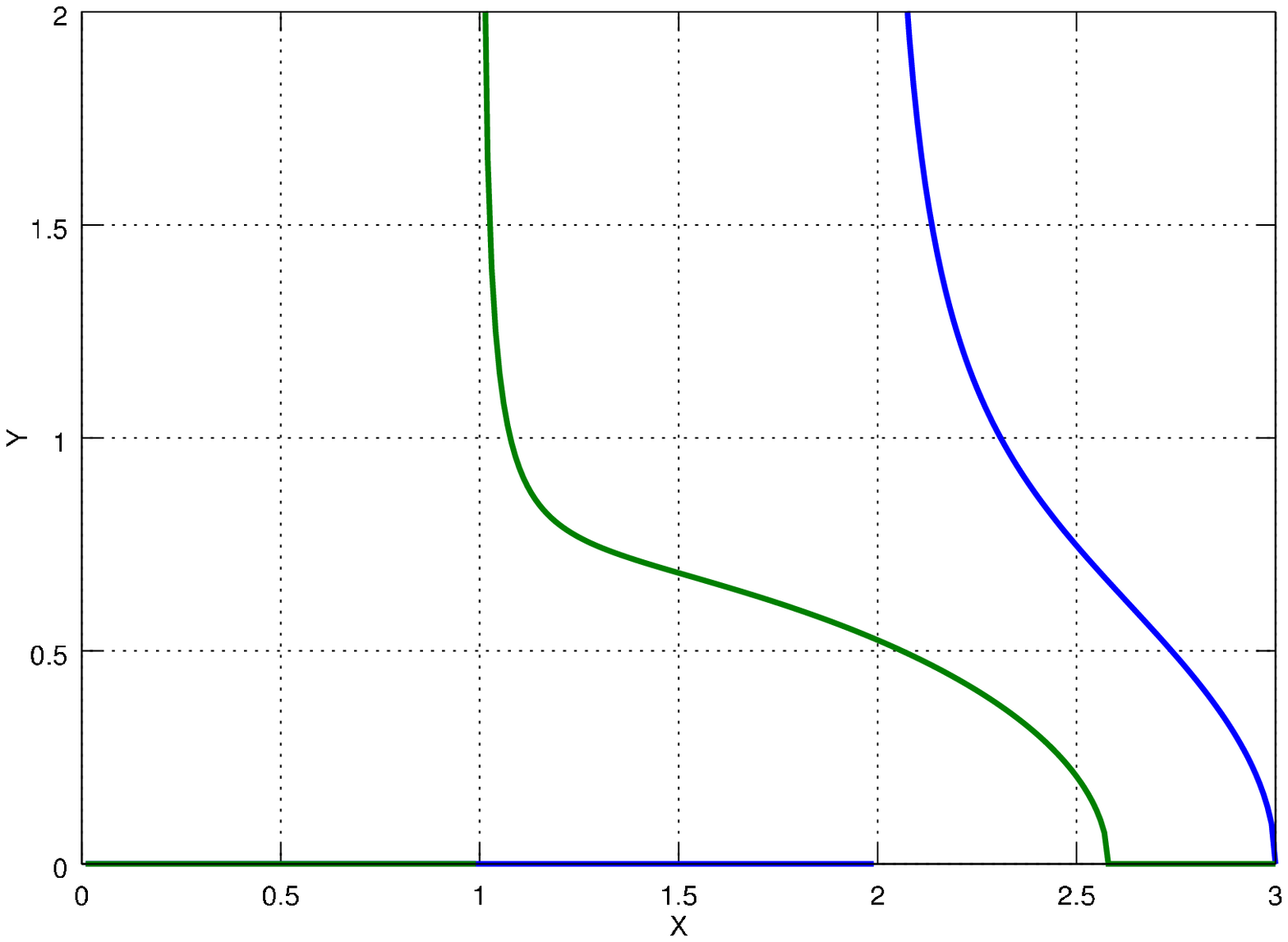}}
\caption{\bf $\alpha=2$.}
\end{figure}

\begin{flushleft}\textcolor{blue}{\rule{1cm}{1pt}} for  $a=2$, $b=3$, $\sigma=2$. Density of eigenvalues in $\Sigma_\sigma=[2,+\infty[$.\\
\textcolor{OliveGreen}{\rule{1cm}{1pt}} for  $a=1$, $b\approx 2.58$, $\sigma=1$.  Density of eigenvalues $\Sigma_\sigma=[1,+\infty[$.\\
\end{flushleft}

\begin{figure}[h]
\centering\scalebox{0.6}{\includegraphics[width=17cm, height=10cm]{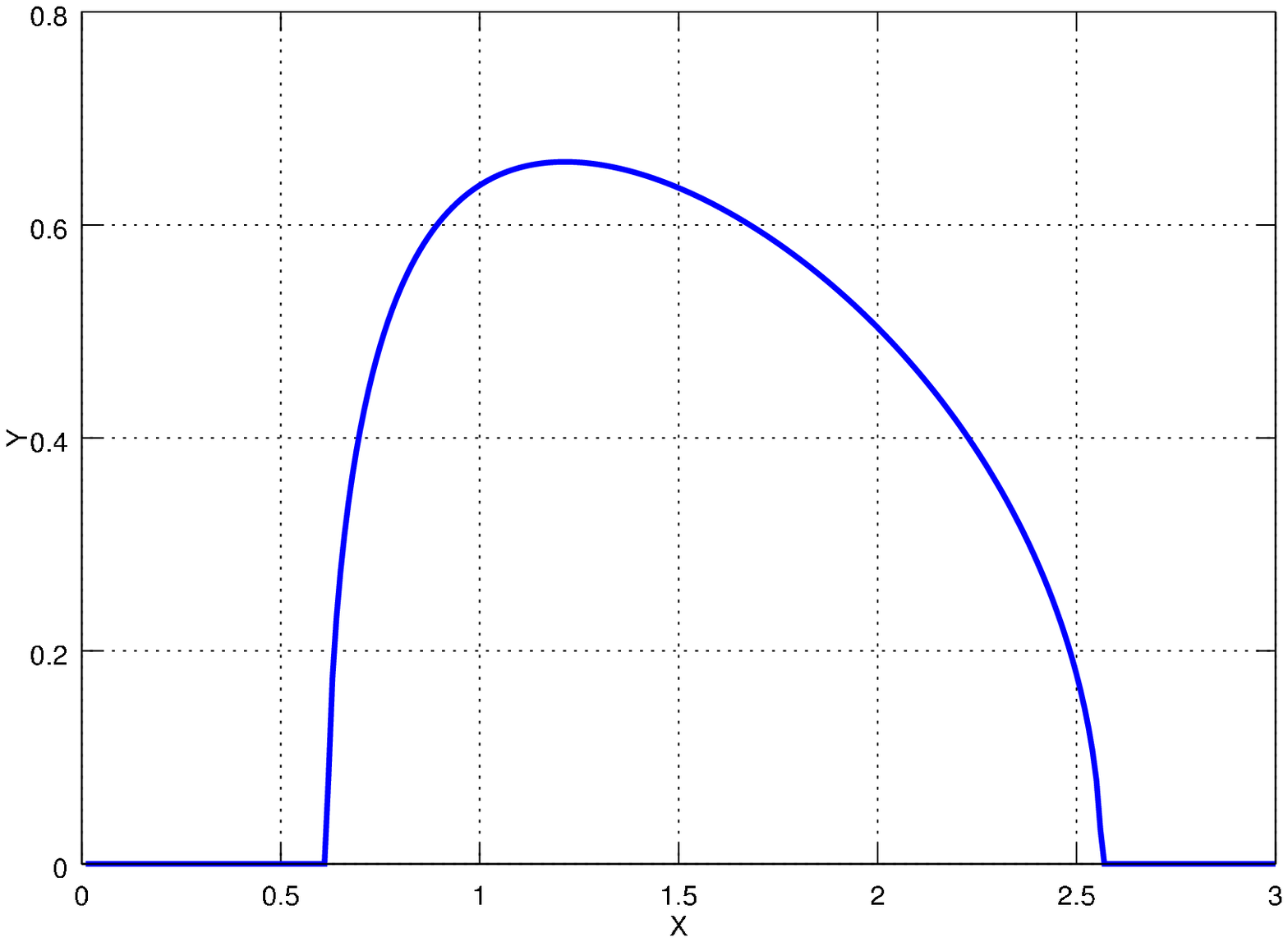}}
\caption{\bf $\alpha=2$, $a=a_c\approx 0.618$, $b=b_c\approx 2.562$}
Density of eigenvalues in $\Sigma_\sigma=]\sigma,+\infty[$, for $0\leq\sigma\leq a_c$
\end{figure}

\begin{figure}[h]
\centering\scalebox{0.5}{\includegraphics[width=17cm, height=10cm]{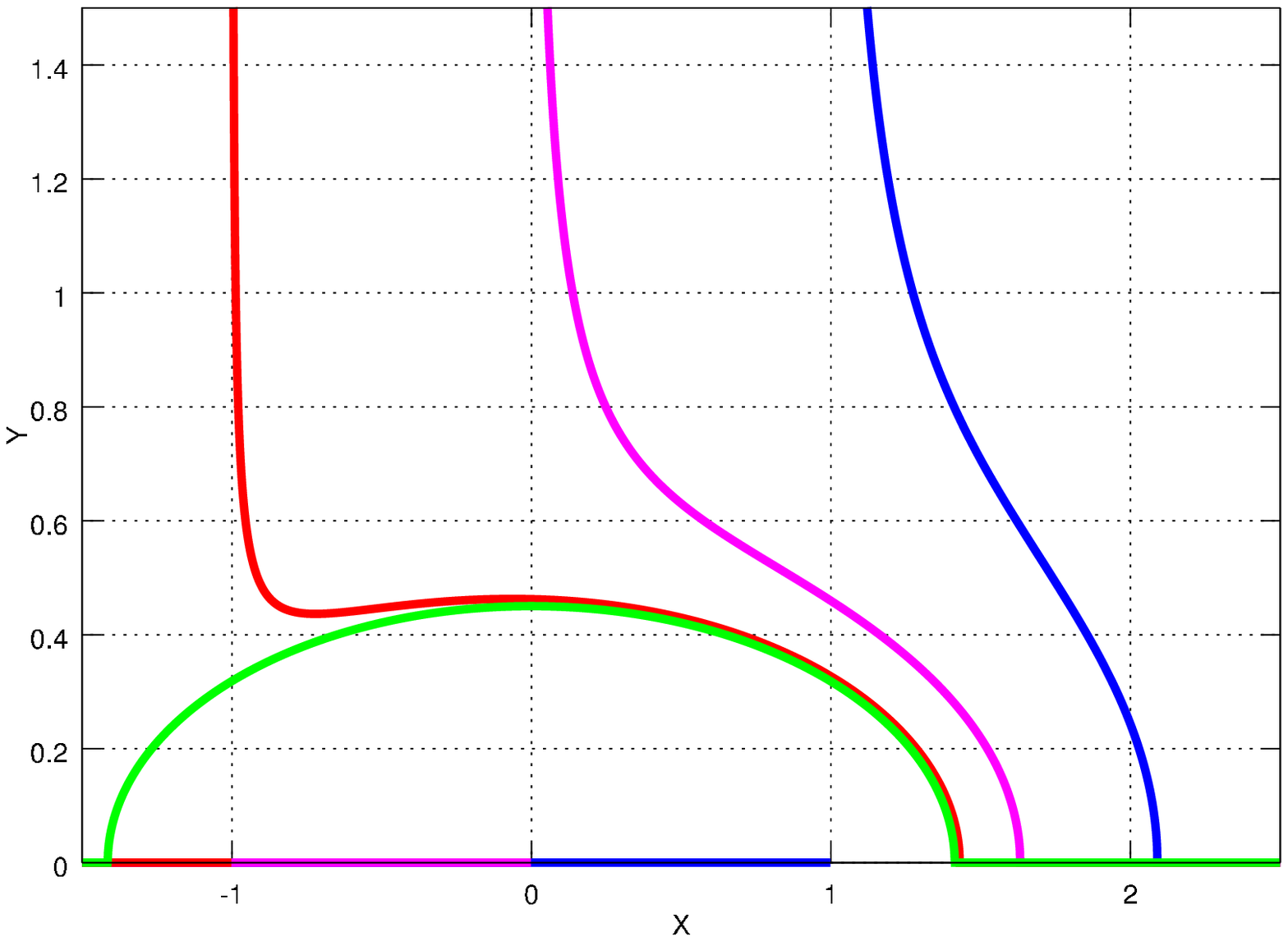}}
\caption{\bf $\alpha=0$}
\end{figure}
\newpage
\begin{flushleft}
\textcolor{blue}{\rule{1cm}{1pt}} for  $a=1$, $b\approx 2.09$, $\sigma=1$, density of eigenvalues in $\Sigma_\sigma=[1,+\infty[$.\\
\textcolor{magenta}{\rule{1cm}{1pt}} for  $a=0$, $b=\frac23\sqrt6\approx 1.632$, $\sigma=0$, density of eigenvalues in $\Sigma_\sigma=[0,+\infty[$.\\
\textcolor{red}{\rule{1cm}{1pt}}for  $a=-1$, $b=1.43$, $\sigma=-1$, density of eigenvalues in $\Sigma_\sigma=[-1,+\infty[$.\\
\textcolor{green}{\rule{1cm}{1pt}}for  $a=-\sqrt 2$, $b=\sqrt 2$, $\sigma\in[-\infty, -\sqrt 2[$, density of eigenvalues in $\Sigma_\sigma=[\sigma,+\infty[$.\\
\end{flushleft}

\begin{figure}[h]
\centering\scalebox{0.5}{\includegraphics[width=17cm, height=14cm]{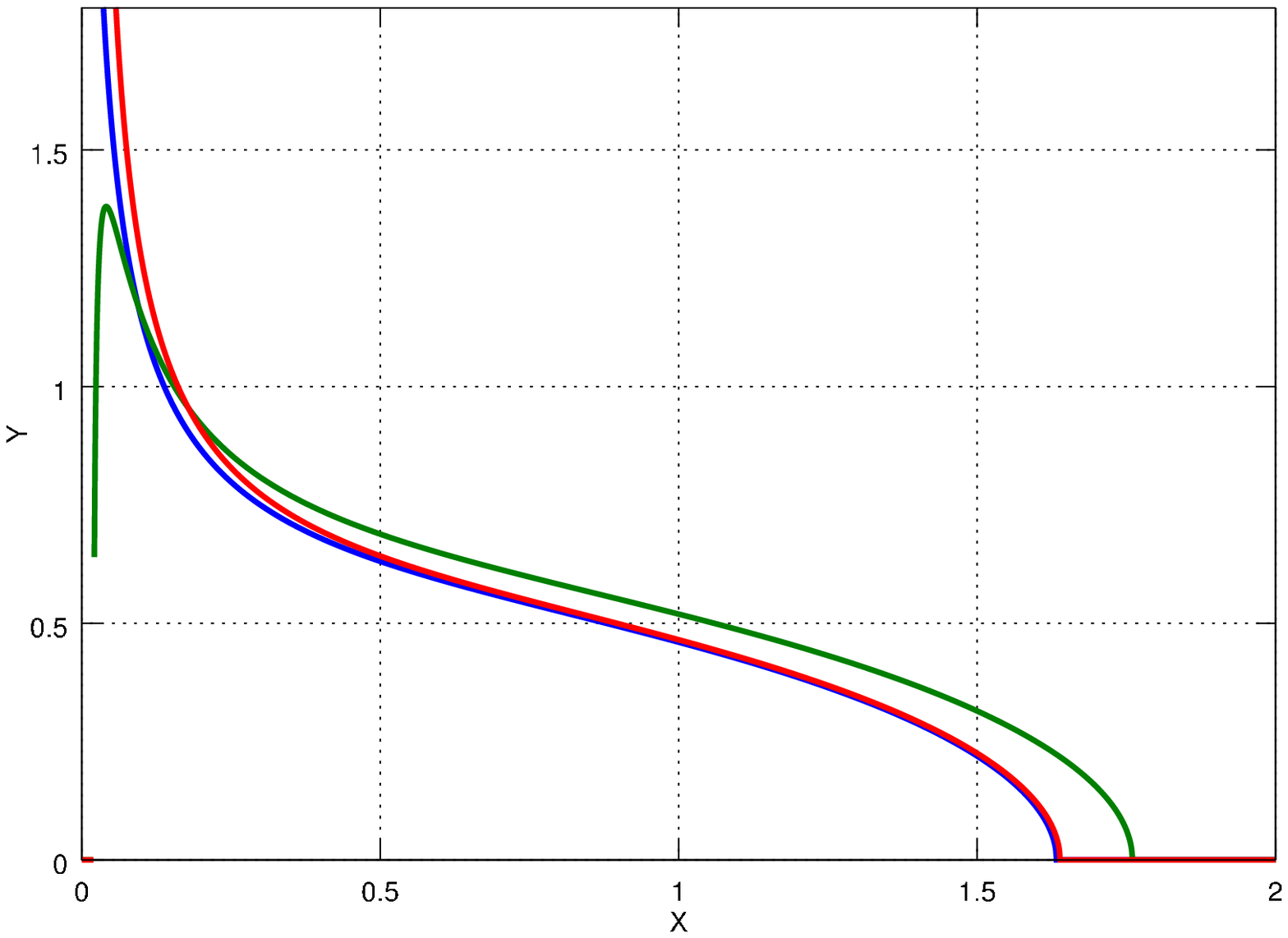}}
\caption{\bf Small value of $\alpha, a$, $\sigma=0$}
Density of eigenvalue in $]\sigma,+\infty[$
\end{figure}\newpage
\begin{flushleft}
\textcolor{OliveGreen}{\rule{1cm}{1pt}} for $2\alpha=0.333.$, $a=a_c\approx 0.02$, $b=b_c\approx 1.7601$ .\\
\textcolor{red}{\rule{1cm}{1pt}} for $\alpha=0$, $a=0.02$, $b\approx 1.639$. The Dean-Majumdar density \\
\textcolor{blue}{\rule{1cm}{1pt}} for $\alpha=0$, $a=0$, $b=\frac23\sqrt6\approx 1.632$ The Dean-Majumdar density.\\
\end{flushleft}

It can be seen that the discontinuity of the density $f_{\alpha, \sigma}$ at $a_c\approx 0$, begin to appear when $a$ and $\alpha$ are near $0$.\\

\begin{center}{\bf Appendix}\end{center}
 {\bf Proof of lemma 3.4.}---\\
 {\bf Step 1:} Substituting $\displaystyle x=\frac{2}{b-a}(t-\frac{a+b}{2})$, it yields
 $$\int_a^b\sqrt{(t-a)(b-t)}\frac{dt}{t}=\frac14(b-a)^2\int_{-1}^1\frac{\sqrt{1-x^2}}{\frac{a+b}{2}+\frac{b-a}{2}x}dx.$$
 now use the change $x=\cos\theta$, $-\pi\leq\theta\leq\pi$ and $\displaystyle r=\frac{b-a}{b+a}\in]0,1]$.
 It yields $$\int_a^b\sqrt{(t-a)(b-t)}\frac{dt}{t}=\frac{r(b-a)}{4}\int_{-\pi}^\pi\frac{\sin^2\theta}{1+
 r\cos\theta}d\theta.$$
 Moreover $$\sin^2\theta=\frac{1}{r^2}\Big(1-r^2\cos^2\theta+r^2-1\Big),$$
hence
 $$\frac{\sin^2\theta}{1+r\cos\theta}=\frac1{r^2}\Big(1-r\cos\theta\Big)+\frac{r^2-1}{r^2}\frac1{1+r\cos\theta}.$$
 Since
 $$1+r\cos\theta=(r+1)\cos^2\frac{\theta}{2}-(r-1)\sin^2\frac{\theta}{2},$$
 hence
 $$\frac{\sin^2\theta}{1+r\cos\theta}=\frac1{r^2}\Big(1-r\cos\theta\Big)+\frac{r-1}{r^2}
 \frac{1+\tan^2\frac{\theta}{2}}{1+\frac{1-r}{r+1}\tan^2\frac{\theta}{2}}.$$
 One can remark that the antiderivative of $\displaystyle\arctan\Big(\sqrt{\frac{1-r}{1+r}}\tan\frac{\theta}{2}\Big)$ is the function $\displaystyle\frac12\sqrt{\frac{1-r}{1+r}}\frac{1+\tan^2\frac{\theta}{2}}{1+\frac{1-r}{1+r}\tan^2\frac{\theta}{2}}.$
 One gets
 $$\int_{-\pi}^\pi\frac{\sin^2\theta}{1+r\cos\theta}d\theta=\frac{2\pi}{r^2}-\frac{2\pi}{r^2}\sqrt{1-r^2}.$$
 Substituting $\displaystyle r=\frac{b-a}{b+a}$, it yields
 $$\int_a^b\sqrt{(t-a)(b-t)}\frac{dt}{t}=\frac{\pi}{2}(\sqrt b-\sqrt a)^2.$$
{\bf Step 2:} For the second step in lemma, it follows from this remark
$$U^\mu(x)=\int_0^{+\infty}\log\frac{1}{|x-t|}\mu(dt)=-\mu\big([0,+\infty[\big)\log(x)+\int_0^{+\infty}\log\frac{1}{|1-\frac{t}{x}|}\mu(dt),$$
and we use first step and dominated convergence theorem.\\
{\bf Step 3:} We will evaluate the integral $\displaystyle \psi(x)=\int_b^x\sqrt{(u-b)(u-a)}\frac{du}{u}$. To do this, take the change of variable $u=\frac{b-a}{2}t+\frac{b+a}{2}$ and $r=\frac{b-a}{b+a}\in]0,1]$, by straightforward computation one gets
$$\psi(x)=\frac{(b-a)^2}{4(b+a)}\int_1^X\frac{\sqrt {t^2-1}}{1+rt}dt,$$ where
$X=\frac{b-a}{2}x+\frac{b+a}{2}.$
Since
$$\begin{aligned}r^2\frac{\sqrt {t^2-1}}{1+rt}&=\frac{1-r^2}{(1+rt)\sqrt{t^2-1}}+\frac{rt-1}{\sqrt{t^2-1}}\\
&=\sqrt{1-r^2}\frac{\sqrt{1-r^2}-r\sqrt{t^2-1}}{(1+rt)\sqrt{t^2-1}}+\frac{rt-1}{\sqrt{t^2-1}}+\frac{r\sqrt{1-r^2}}{1+rt}
,\end{aligned}$$
Moreover \begin{equation}\varphi_1(x):=\frac{rt-1}{\sqrt{t^2-1}}=\frac{rt}{\sqrt{t^2-1}}-\frac{1+\frac{t}{\sqrt{t^2-1}}}{t+\sqrt{t^2-1}},\end{equation}
and $$(1+rt)\Big(\sqrt{1-r^2}-r\sqrt{t^2-1}\Big)=-\Big(\sqrt{t^2-1}-t\sqrt{1-r^2}\Big)\Big(t+r+\sqrt{1-r^2}\sqrt{t^2-1}\Big),$$
hence
$$\frac{\sqrt{1-r^2}-r\sqrt{t^2-1}}{(1+rt)\sqrt{t^2-1}}=-\frac{\Big(\sqrt{t^2-1}-t\sqrt{1-r^2}\Big)\Big(t+r+\sqrt{1-r^2}
\sqrt{t^2-1}\Big)}{(1+rt)^2\sqrt{t^2-1}},$$
which can be written as
$$\frac{\sqrt{1-r^2}-r\sqrt{t^2-1}}{(1+rt)\sqrt{t^2-1}}=-\frac{\sqrt{t^2-1}-t\sqrt{1-r^2}}{\sqrt{t^2-1}\Big(t+r-\sqrt{1-r^2}\sqrt{t^2-1}\Big)},$$
hence \begin{equation}\varphi_2(x):=\sqrt{1-r^2}\frac{\sqrt{1-r^2}-r\sqrt{t^2-1}}{(1+rt)\sqrt{t^2-1}}=-\frac{1-\sqrt{1-r^2}\frac{t}{\sqrt{t^2-1}}}
{t+r-\sqrt{1-r^2}\sqrt{t^2-1}}.\end{equation}
Thus $$\frac{\sqrt {t^2-1}}{1+rt}=\frac1{r^2}\Big(\varphi_1(t)+\varphi_2(t)+\frac{r\sqrt{1-r^2}}{1+rt}\Big),$$
and all the functions can be integrate explicitly and one gets
$$\begin{aligned}\psi(x)&=\frac1{r^2}\Big(r\sqrt{X^2-1}-\sqrt{1-r^2}\log\big(X+\sqrt{X^2-1}\big)-\sqrt{1-r^2}\log\big(X+r-\sqrt{1-r^2}\sqrt{X^2-1}\big)\\&+\sqrt{1-r^2}\log(1+rX)
\Big).\end{aligned}$$
By performing the change of variable in the reverse order $X=\frac{b-a}{2}x+\frac{b+a}{2}$, and expanding near $+\infty$, $\psi(x)$, one gets
$$\psi(x)=x-\frac12(b+a)\log x+A+o(\frac1x),$$
where $$\displaystyle A=-\frac{a+b}{2}\big(1-\log\frac{b-a}{4}\big)+\sqrt{ab}\log\frac{\sqrt b-\sqrt a}{\sqrt b+\sqrt a}.$$
Using step 2, one gets $C_\mu=-A$. This complete the proof.\\
\begin{lemma} For $0<a<b$ and as $x$ goes to $+\infty$, we have
\begin{description}
\item [(1)] $$\int_b^x\sqrt{(b-u)(u-a)}\frac{du}{u}=x-\frac12(b+a)\log x-\frac{a+b}{2}\big(1-\log\frac{a+b}{4}\big)+\sqrt{ab}\log\frac{\sqrt b-\sqrt a}{\sqrt b+\sqrt a}+o(1),$$
\item[(2)]$$\int_b^x\sqrt{(t-a)(t-b)}dt=\frac12x^2-\frac{a+b}{2}x-\frac{(b-a)^2}{8}\log x+\frac1{16}\Big(a^2+6ab+b^2+2(b-a)^2\log\frac{b-a}{4}\Big)+o(1).$$
\item[(3)]$$\int_b^x\sqrt{\frac{t-b}{t-b}}dt=x-\frac{b-a} 2\log x-\frac{a+b}{2}+\frac{b-a}{2} \log\frac{b-a}{4}+o(1).$$
\end{description}
 \end{lemma}
{\bf Proof}.--- \\
{\bf Step(1):} This step is proved in the previous lemma.\\
{\bf Step (2):}Take the change of variable $u=\frac{t-a}{b-a}$, hence
$$\int_b^x\sqrt{(t-a)(t-b)}dt=(b-a)^2\int_1^X\sqrt{u(u-1)}du,$$ where $X=\frac{x-a}{b-a}$.
Since for all $u>1$, $$\sqrt{u(u-1)}=\frac{1}{4} \Big(-\frac{1}{2 \sqrt{-1+u} \sqrt{u}}+2 \sqrt{(-1+u) u}+\frac{(-1+2 u)^2}{2 \sqrt{(-1+u) u}}\Big),$$
moreover $$\int \frac{1}{2 \sqrt{-1+u} \sqrt{u}}du={\rm arcsinh}(\sqrt{u-1}),$$
 and $$\int\Big(2 \sqrt{(u-1) u}+\frac{(-1+2 u)^2}{2 \sqrt{(u-1) u}}\Big)du=(2u-1)\sqrt{u(u-1)}.$$
 Hence
 $$\int_1^X\sqrt{u(u-1)}du=\frac14\Big((2X-1)\sqrt{X(X-1)}-{\rm arcsinh}(\sqrt{X-1})\Big).$$
 Substituting the expression of $X$ and as $x$ goes to $+\infty$ one gets
 $$\int_b^x\sqrt{(t-a)(t-b)}dt=\frac12x^2-\frac{a+b}{2}x-\frac{\theta^2}{8}\log x+\frac1{16}\Big(a^2+6ab+b^2+2(b-a)^2\log\frac{b-a}{4}\Big)+o(1).$$
{\bf Step (3):} As in the previous by performing the change of variable $u=\frac{t-a}{b-a}$ and $X=\frac{x-a}{b-a}$, one obtains
$$\int_b^x\sqrt{\frac{t-b}{t-a}}dt=(b-a)\int_1^X\sqrt{\frac{u-1}{u}}du.$$
Since $$\sqrt{\frac{u-1}{u}}=-\frac{\frac{1}{2 \sqrt{-1+u}}+\frac{1}{2 \sqrt{u}}}{\sqrt{-1+u}+\sqrt{u}}+\frac{-1+2 u}{2 \sqrt{(-1+u) u}},$$
hence $$\int_1^X\sqrt{\frac{u-1}{u}}du=\sqrt{X(X-1) }-\log\big(\sqrt{X-1}+\sqrt{X}\big).$$
Substituting the change of $X$ and when $x$ goes to $+\infty$, one gets the desired result.\\
{\bf Approximation Density.} Recall that for $\lim_{n\to+\infty}\frac{\mu_n}{n}=\alpha$, the density of positive eigenvalues is $\displaystyle f_\alpha(x)$ is given in theorem 2.4.

Let $$f_n(x)=\frac1{\sqrt n}\sum_{k=0}^{n-1}\varphi^{\mu_n}_k(\sqrt nx)^2,$$
where $\displaystyle\varphi^{\mu_n }_k(x)=\frac1{\sqrt {d_{k,n}}}H^{\mu_n}_k(x)x^{\mu_n}e^{-\frac{x^2}{2}},$  and $H_k^{\mu_n}$ is the truncated orthogonal Hermite polynomial on the positive real axis, which satisfies $$\int_0^{+\infty}H_k^{\mu_n}(x)H_m^{\mu_n}(x)x^{2\mu_n}e^{-x^2}dx=0,\quad {\rm for}\; m \neq k,$$
and $$\int_0^{+\infty}(H_k^{\mu_n}(x))^2x^{2\mu_n}e^{-x^2}dx=d_{k, n}.$$
It has been proved in theorem 2.4, that as n go to $+\infty$, the density $f_n$ converge tightly to the density $f_\alpha$ where $\displaystyle\alpha=\lim_{n\to\infty}\frac{\mu_n}n$.

{\bf First case:} $n=7, \mu_7=0$, hence $\alpha=0$.
$$\displaystyle f_0(x)=\frac1{2\pi}\sqrt{\frac{\frac23\sqrt 6-x}{x}}\big(2x+\frac23\sqrt 6\big), {\rm with\;\; support}\;\; [0,\frac23\sqrt 6]$$.

\begin{figure}[h]
\centering\scalebox{0.5}{\includegraphics[width=17cm, height=14cm]{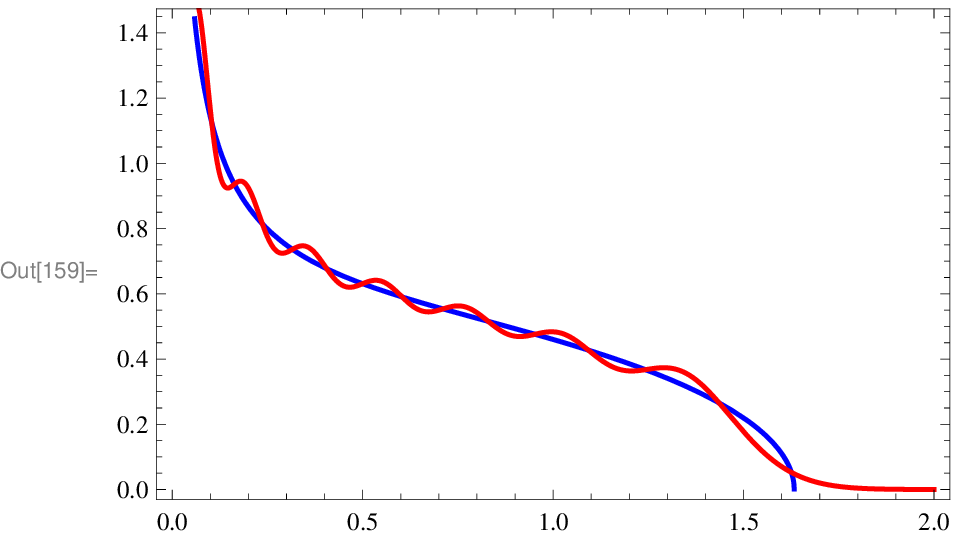}}

\end{figure}
\begin{flushleft}
\textcolor{blue}{\rule{1cm}{1pt}} Exact density $f_0$ of positive eigenvalues .\\
\textcolor{red}{\rule{1cm}{1pt}} Approximative density $f_7$ of positive eigenvalues.\\
\end{flushleft}
\newpage
{\bf Second case:} $n=5, \mu_5=\frac52$, hence $\alpha=\frac12$.
$$\displaystyle f_{\frac12}(x)=\frac1{\pi}\sqrt{(1.9-x)(x-0.1}\big(1+\frac1{2x\sqrt{0.19}}\big), {\rm with\;\; support}\;\; [0.1,1.9].$$
\begin{figure}[h]
\centering\scalebox{0.5}{\includegraphics[width=17cm, height=14cm]{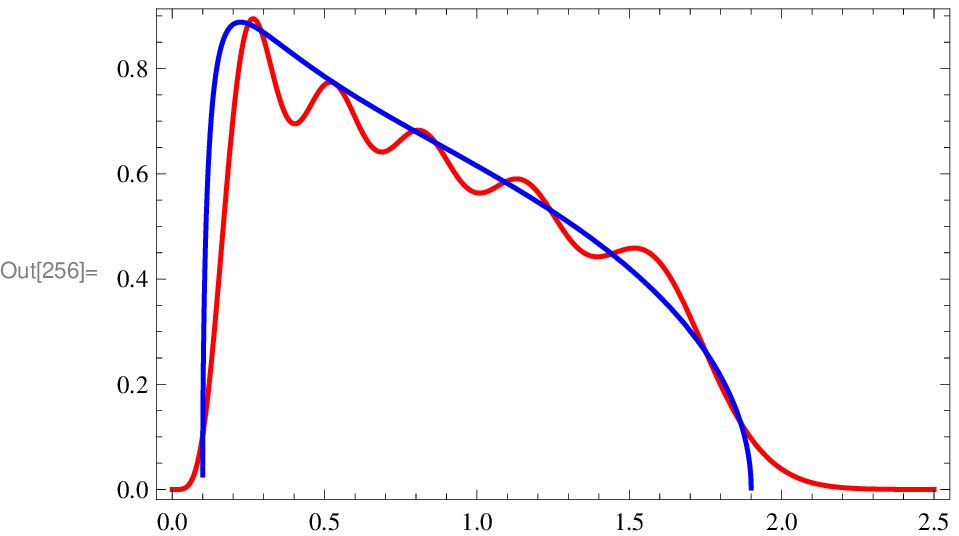}}

\end{figure}
\begin{flushleft}
\textcolor{blue}{\rule{1cm}{1pt}} Exact density $f_{\frac12}$ of positive eigenvalues .\\
\textcolor{red}{\rule{1cm}{1pt}} Approximative density $f_5$ of positive eigenvalues.\\
\end{flushleft}

\begin{center}{\bf Acknowledgments.}
\end{center}
My sincere thanks go to Jacques Faraut for his comments on this manuscript and his important remarks.

Address:  College of Applied Sciences
   Umm Al-Qura University
  P.O Box  (715), Makkah,
  Saudi Arabia.\\
E-mail: bouali25@laposte.net \& mabouali@uqu.edu.sa

\end{document}